\documentclass[a4paper,oneside,11pt]{article}

\usepackage{amsfonts,color}

\newtheorem{dfn}{Definition}[section]
\newtheorem{tw}[dfn]{Theorem}
\newtheorem{prop}[dfn]{Proposition}
\newtheorem{rem}[dfn]{Remark}
\newtheorem{ex}[dfn]{Example}
\newtheorem{lem}[dfn]{Lemma}

\newtheorem{cor}[dfn]{Corollary}
\usepackage{amssymb} 
\usepackage{amsmath}
\numberwithin{equation}{section}

 \global\long\def\sbr#1{\left[ #1\right] }
 \global\long\def\cbr#1{\left\{  #1\right\}  }
 
 \global\long\def\rbr#1{\left(#1\right)}

 \global\long\def\R{\mathbb{R}}
\global\long\def\S{\mathbb{S}}

 \global\long\def\dd#1{\textnormal{d}#1}

 \global\long\def\ra{\rightarrow}
 
 \global\long\def\ns{\infty}

\usepackage{eurosym}

\usepackage{anysize}
\usepackage{array}
\usepackage{enumerate}

\author{Micha\l \ Barski  \\ \small  Faculty of Mathematics, Warsaw University, Poland\\
 \small{\it m.barski@mimuw.edu.pl} \bigskip \\
\\
Rafa\l \ \L ochowski
\\ \small 
Department of Mathematics and Mathematical Economics,\\ \small Warsaw School of Economics, Poland\\ \small{\it rlocho@sgh.waw.pl}}

\title{\bf On the reducibility of affine models with dependent L\'evy factors}

\begin{document}

\maketitle

\begin{abstract}
The paper is devoted to the study of the short rate equation of the form
\begin{gather*}
\dd R(t)=F(R(t))\dd t+\sum_{i=1}^{d}G_i(R(t-))\dd Z_i(t), \quad R(0)=x\geq 0,\quad t>0,
\end{gather*}
with deterministic functions $F,G_1,...,G_d$ and a multivariate L\'evy process $Z=(Z_1,...,Z_d)$ with possibly dependent coordinates. The equation is supposed to have a nonnegative solution which generates an affine term structure model. The L\'evy measure $\nu$ of $Z$ is assumed to admit a spherical decomposition based on the representation $\mathbb{R}^d=\S^{d-1}\times (0,+\infty)$, where $\S^{d-1}$ stands for the unit sphere. Then $\nu(\dd y)=\lambda(\dd \xi)\times \gamma_{\xi}(\dd r)$, where $\lambda$ is a measure on $\S^{d-1}$ and $\gamma_{\xi}$ on $(0,+\infty)$. Under some assumptions on spherical decomposition, a precise form of the generator of $R$ is determined and it is shown that the resulted term structure model is identical to that generated by the equation
$$
\dd R(t)=(a R(t)+b)\dd t+C\cdot (R(t-))^{1/\alpha} dZ^{\alpha}(t),  \quad R(0)=x,
$$
with some constants $a,b,C$ and a one dimensional $\alpha$-stable L\'evy process $Z^{\alpha}$, where $\alpha\in(1,2)$. 
The case when $\nu$ has a density is considered as a special case.

The paper generalizes the classical results on the Cox-Ingersoll-Ross (CIR) model, \cite{CIR}, as well as on its extended version from \cite{BarskiZabczykCIR} and \cite{BarskiZabczyk} where $Z$ is a one-dimensional L\'evy process. It is the starting point for the 
classification in the spirit of \cite{DaiSingleton} and \cite{BarskiLochowski} for the affine models with dependent L\'evy processes.

\end{abstract}


\section{Introduction}

Let us consider a bond market with a family of stochastic processes describing zero coupon bond prices 
$
P(t,T), t\in[0,T], 
$
parametrized by the maturity time $T>0$ and the short rate process
$
R(t), t\geq 0.
$
The processes are defined on a probability space $(\Omega,\mathcal{F},\mathbb{P})$ with a filtration $(\mathcal{F}_t), t\geq 0$.
The bond with maturity $T$ pays to its owner at time $T$ a nominal value  assumed here to be $1$, i.e. $P(T,T)=1$.
The discounted value of  $1$ paid at time $t>0$ equals 
$
D(t)=e^{-\int_{0}^{t}R(s)ds}.
$
The short rate process $R$ is supposed to satisfy, for each $T>0$, the condition
\begin{gather}\label{wwo exp affine}
\mathbb{E}[e^{-\int_{t}^{T}R(s)ds}\mid \mathcal{F}_t]=e^{-A(T-t)-B(T-t)R(t)}, \qquad t\in[0,T], 
\end{gather}
with some deterministic functions $A(\cdot), B(\cdot)$. Interpreting $\mathbb{P}$ as a risk neutral measure, one recognizes in the left side of \eqref{wwo exp affine} the price at time $t$ of the bond with maturity $T$, that is $P(t,T)$. Thus \eqref{wwo exp affine} means that 
the short rate $R$ {\it generates} an {\it affine term structure}.

The concept of modelling bond prices in the affine fashion was introduced by Filipovi\'c \cite{FilipovicATS} and Duffie, Filipovi\'c and 
Schachermeyer \cite{DuffieFilipovicSchachermeyer}. It was motivated by the results of Kawazu and Watanabe \cite{KawazuWatanabe}
on continuous state branching processes with immigration. Further developments on regularity of affine processes 
are due to Cuchiero, Filipovi\'c and Teichmann \cite{CuchieroFilipovicTeichmann} and Cuchiero and Teichmann
\cite{CuchieroTeichmann}. The aforementioned results are settled in the general Markovian setting and the description of affine processes 
is given in the form of their generators. A class of particular interest are short rates given by stochastic equations. 
We call {\it a generating equation} an equation with solution which generates an affine model. The precursors of generating equations are two classical equations due to
 Cox-Ingersoll-Ross (CIR) \cite{CIR}  
\begin{gather}\label{CIR}
\dd R(t)=(aR(t) +b)\dd t+c\sqrt{R(t)}\dd W(t), 
\end{gather}
with $a\in\mathbb{R}, b\geq 0,c>0$
and due to Vasi\v cek \cite{Vasicek} 
\begin{gather}\label{Vasicek}
\dd R(t)=(aR(t) +b)\dd t+c \ \dd W(t), 
\end{gather}
with $a,b,c\in\mathbb{R}$, both driven by a one dimensional Wiener process $W$. To make the behaviour of the short rate process more realistic and to improve the accuracy of calibration to market data more involved equations are considered in the literature. 
Dai and Singleton \cite{DaiSingleton} consider factorial models perturbed by correlated Wiener processes and examine the influence of the correlation structure on the resulting affine model. It is shown in \cite{DaiSingleton}, in particular, that different equations may generate identical affine models.
Replacing $W$ by a L\'evy process it was shown in Barski and Zabczyk
\cite{BarskiZabczykCIR} that the generalization of \eqref{CIR} must be of the form
\begin{gather}\label{rownanie na R stabilny d=1}
\dd R(t)=(a R(t)+b)\dd t+C\cdot (R(t-))^{1/\alpha} dZ^{\alpha}(t), 
\end{gather}
with $a\in\mathbb{R}, b\geq 0,C>0$, where $Z^{\alpha}$ is an $\alpha$-stable process with index $\alpha\in(1,2]$. It was also shown in \cite{BarskiZabczykCIR} that the counterpart of \eqref{Vasicek} in the L\'evy setting allows to preserve the positivity of $R$, which clearly is not preserved in \eqref{Vasicek}. Models driven by a multivariate L\'evy process with independent coordinates appear, among others, in Duffie and G\^arleanu \cite{DuffieGarleanu},  Barndorff-Nielsen and Shephard \cite{Bandorff-NielsenShepard}, Jiao, Ma and Scotti \cite{JiaoMaScotti}, Barski and \L ochowski \cite{BarskiLochowski}. Similarly as in \cite{DaiSingleton}, it is noticed in \cite{BarskiLochowski} that different equations may generate identical affine models. This fact motivated a classification of all generating equations into several classes which are representable by the so called canonical equations having tractable forms.
 The case when the coordinates of the multivariate L\'evy process are dependent is an unexplored field entered by this paper.

We consider an equation of the form
\begin{gather}\label{main equation}
\dd R(t)=F(R(t))\dd t+\sum_{i=1}^{d}G_i(R(t-))\dd Z_i(t), \quad R(0)=x\geq 0,\quad t>0,
\end{gather}
where $F,G:=(G_1,...,G_d)$ are deterministic functions and $Z=(Z_1,...,Z_d)$ is a multivariate L\'evy process. As the coordinates of $Z$ may be dependent, the L\'evy measure $\nu$ of $Z$ is not necessarily concentrated on axes and it is assumed to admit  
a {\it spherical decomposition} which means that 
\begin{equation} \label{spherical}
\nu(A)=\int_{\S^{d-1}}\int_{0}^{+\infty}\mathbf{1}_{A}(r\xi)\gamma_{\xi}(\dd r) \ \lambda(\dd \xi), \quad A\in\mathcal{B}(\mathbb{R}^d).
\end{equation}
 Here $\S^{d-1}:= \{x\in\mathbb{R}^d:\mid x\mid=1\}$ is a unit sphere in $\mathbb{R}^d$, $\lambda$ is a finite measure on $\S^{d-1}$ called a {\it spherical part} of $\nu$ and $\{{\gamma_{\xi}};{\xi\in\S^{d-1}}\}$ is a family of measures on $(0,+\infty)$ called {\it radial measures} of $\nu$. The spherical decomposition may seem a very particular feature of $\nu$
at first glance, but  let us recall a well known example of an $\alpha$-stable process in $\mathbb{R}^d$ with index $\alpha\in(1,2)$. Its jump measure admits a spherical decomposition with identical radial measures given by the density
\begin{gather}\label{radialna stabilna}
\gamma_{\xi}(\dd r)=\gamma(\dd r):=\frac{1}{r^{1+\alpha}}\dd r, \quad r>0,  \ \xi\in\S^{d-1},
\end{gather}
and with any finite spherical measure $\lambda$ on $\S^{d-1}$. In fact, the spherical decomposition can be determined, for instance, when $\nu$ has a density with respect to the Lebesgue measure. In this case the radial measures are of the form
$$
\gamma_{\xi}(\dd r)=g(r\xi) r^{d-1}\sqrt{1-\xi_1^2}\cdot \sqrt{1-(\xi_1^2+\xi_2^2)}\cdot...\cdot
\sqrt{1-(\xi_1^2+...+\xi_{d-2}^2)} \ \dd r.
$$
for $\xi=(\xi_1,...,\xi_d)\in\S^{d-1}, r>0$,	 where $g$ stands for the density of $\nu$. The spherical part $\lambda$ is then the image of the Lebesgue measure by the polar transformation, see Section \ref{sec o Laplacach dla gestosci wiki} for details. We examine the question which affine models can be generated by equation \eqref{main equation} 
when \eqref{spherical} holds. In Example  \ref{ex1} we show that 
if \eqref{main equation} is a generating equation and $Z$ is an $\mathbb{R}^d$-valued $\alpha$-stable process then the resulted affine model is identical with the model generated by \eqref{rownanie na R stabilny d=1}. 
This means that \eqref{rownanie na R stabilny d=1} can replace the initial equation, which may be of a complicated form, preserving the bond prices unchanged. In this case we call the initial equation to have the {\it reducibility property}. This extends the observations from \cite{DaiSingleton} and \cite{BarskiLochowski} to the case with dependent noise coordinates. The main result of this paper, Theorem \ref{tw_sferyczne}, provides conditions for \eqref{main equation} driven by $Z$ with L\'evy measure of the form \eqref{spherical} to have  the reducibility property. We prove that if $G$ is a continuous function for which the limit 
$$
\lim_{x\rightarrow 0^{+}}\frac{G(x)}{\left| G(x)\right|}
$$
exists and the Laplace exponents related to radial measures 
$$
J_{\gamma_{\xi}}(b):=\int_{0}^{+\infty}(e^{-br}-1+br)\gamma_{\xi}(\dd r), \quad b\geq 0, \quad \xi\in\S^{d-1},
$$
satisfy the condition
\begin{gather}\label{warunek na Laplasy}
\sup_{\xi\in\text{supp}(\lambda)}J_{\gamma_\xi}(b)\leq K \cdot \inf_{\xi\in\text{supp}(\lambda)}J_{\gamma_\xi}(b), \quad b\geq 0,
\end{gather}
with some $1\leq K<+\infty$, then any equation with such $G$ and $Z$ can generate only the same affine model as that one generated by \eqref{rownanie na R stabilny d=1} with some $\alpha\in(1,2)$. In particular, this takes place, for instance, when $Z$ is $\mathbb{R}^2$-valued and its jump measure has a density $g$ such that the functions
$$
\underline{g}(r) := \inf_{\left|x\right| = r}g(x), \qquad \bar{g}(r) := \sup_{\left|x\right| = r}g(x),  \quad r\geq 0,
$$ 
satisfy the integrability conditions
$$
0<\int_{0}^{+\infty}\min\{r^2,r^3\}\underline{g}(r)\dd r\leq \int_{0}^{+\infty}\min\{r^2,r^3\}\bar{g}(r)\dd r<+\infty,
$$
and 
$$
\lim_{\varepsilon\downarrow 0}\frac{\int_{\varepsilon}^{1}r^2 \bar{g}(r)\dd r}{\int_{\varepsilon}^{1}r^2 \underline{g}(r)\dd r}<+\infty \quad \text{and} \quad \lim_{\varepsilon\downarrow 0}\frac{\int_{1}^{1/\varepsilon}r^3 \bar{g}(r)\dd r}{\int_{1}^{1/\varepsilon}r^3 \underline{g}(r)\dd r}<+\infty,
$$
for details and extensions for the cases $d>2$ see Section \ref{sec o Laplacach dla gestosci wiki}.

The paper is organized as follows. In Section \ref{section Formulation of the problem} we present some basic facts on L\'evy processes and 
the Markovian characterization of generating equations. The main Section \ref{sec Reducibility of equations with multivariate noise} on the reducibility problem is divided into three sections. The first contains the formulation of the main results including Theorem \ref{tw_sferyczne} and illustrative analysis of the case when $\nu$ has a density resulting in Theorem \ref{th reducibility density}. The next section consists of a sequence of auxiliary results and in the final one we present the proof of Theorem \ref{tw_sferyczne}.

\section{Preliminaries}\label{section Formulation of the problem}

\subsection{Basic facts on L\'evy processes}\label{sec basic facts on levy processes}
Let $Z:=(Z_1,Z_2,...,Z_d)$  be a L\'evy process in $\mathbb{R}^d, d\geq1 $ on some probability space $(\Omega,\mathcal{F},\mathbb{P})$ with a  filtration $\{\mathcal{F}_t, t\geq 0\}$. If $Z$ is a martingale, then it admits the following unique representation 
$$
Z(t)=W(t)+X(t), \quad t\geq 0,
$$
where $W$ is a Wiener process in $\mathbb{R}^d$ with a covariance matrix $Q$ and $X$ is the so called {\it jump martingale part}  of $Z$. It is independent of $W$ and can be described in terms of the {\it jump measure} of $Z$ defined by
\begin{gather}\label{jump measure}
\pi(t,A):=\sharp\{s\in[0,t]: \triangle Z(s)\in A\}, \quad t\geq 0,
\end{gather}
where $\triangle Z(s):=Z(s)-Z(s-)$ and $A\subset \mathbb{R}^d$ is a set separated from zero, i.e. $0$ does not belong to the closure of $A$.
With \eqref{jump measure} at hand one defines the {\it L\'evy measure} of $Z$ by
$$
\nu(A):=\mathbb{E}\Big[\pi(1,A)\Big].
$$
Then $X$ can be written as
$$
X(t):=\int_{0}^{t}\int_{\mathbb{R}^d}y \ (\pi(\dd s,\dd y)- \dd s \ \nu(\dd y)), \quad t\geq 0,
$$
and its properties can be formulated in terms of the measure $\nu$. The integrability of $X$ is equivalent 
to the condition
\begin{gather}\label{Levy measure for a martingale}
\int_{\mathbb{R}^d}(\left| y\right|^2\wedge \left|y\right|)\nu(\dd y)<+\infty,
\end{gather}
while the {\it variation} of $X$ is finite if and only if
\begin{gather}\label{Levy measure finite variation}
\int_{\left| y\right|<1} \left| y\right| \nu(\dd y)<+\infty.
\end{gather}
In our notation $\left|\cdot\right|$ stands for the standard norm in $\mathbb{R}^d$ and $\langle\cdot,\cdot\rangle$ for the standard scalar product.

By the independence of $X$ and $W$ we see that, for $\lambda\in\mathbb{R}^d$,
$$
\mathbb{E}\sbr{e^{-\langle \lambda,Z(t)\rangle}}=\mathbb{E}\sbr{e^{-\langle \lambda,W(t)\rangle}}\cdot\mathbb{E}\sbr{e^{-\langle \lambda,X(t)\rangle}},
$$
so the Laplace exponent $J_Z$ of $Z$ defined by 
$$
\mathbb{E}\sbr{e^{-\langle \lambda,Z(t)\rangle}}=e^{tJ_{Z}(\lambda)},
$$
exists at $\lambda$ if and only if $J_X(\lambda)$ finite. The latter property 
is equivalent to the condition
\begin{gather}\label{exp moments}
\int_{\mid y\mid>1}e^{-\langle \lambda,y\rangle}\nu(\dd y)<+\infty.
\end{gather}
If \eqref{exp moments} holds, then 
\begin{equation} \label{Jdef}
J_X(\lambda)=\int_{\mathbb{R}^d}(e^{-\langle\lambda,y\rangle}-1+\langle\lambda,y\rangle)\nu(\dd y),
\end{equation}
and, consequently,
\begin{align} \label{LaplaceZ}\nonumber
J_Z(\lambda)&=J_{W}(\lambda)+J_X(\lambda)\\[1ex]
&={{\frac{1}{2}\langle Q\lambda,\lambda\rangle+\int_{\mathbb{R}^d}(e^{-\langle\lambda,y\rangle}-1+\langle\lambda,y\rangle)\nu(\dd y)}}.
\end{align}
It follows, in particular, that the process $Z$ is uniquely determined by the pair $(Q,\nu)$.

\subsection{Markovian characterization of generating equations}
It was shown in \cite [Theorem 5.3]{FilipovicATS} that the generator of a general positive Markovian short rate generating  an affine model is of the form
\begin{align}\label{generator Filipovica}
\mathcal{A}f(x)=&c x f^{\prime\prime}(x)+(\beta x+\gamma)f^\prime(x)\\[1ex] \nonumber
&+\int_{(0,+\infty)}\Big(f(x+v)-f(x)-f^\prime(x)(1\wedge v)\Big)(m(\dd v)+x\mu(\dd v)), \quad x\geq 0,
\end{align}
for $f\in\mathcal{L}(\Lambda)\cup C_c^2(\mathbb{R}_{+})$, where 
$\mathcal{L}(\Lambda)$ is the linear hull of $\Lambda:=\{f_\lambda:=e^{-\lambda x}, \lambda\in(0,+\infty)\}$
and $C_c^2(\mathbb{R}_{+})$ stands for the set of twice continuously differentiable functions with compact support in $[0,+\infty)$. 
Above $c, \gamma\geq 0$, $\beta\in\mathbb{R}$ and $m(\dd v)$, $\mu(\dd v)$ are nonnegative Borel measures on $(0,+\infty)$ satisfying
\begin{gather}\label{warunki na iary Filipovica}
\int_{(0,+\infty)}(1\wedge v)m(\dd v)+\int_{(0,+\infty)}(1\wedge v^2)\mu(\dd v)<+\infty.
\end{gather}
Moreover, the functions $A(\cdot), B(\cdot)$ in \eqref{wwo exp affine} are uniquely determined by the form of the generator 
\eqref{generator Filipovica}, for details see \cite{FilipovicATS}.

Application of the characterization above to the case when $R$ is given by \eqref{main equation} leads to necessary and sufficient conditions
making  \eqref{main equation}  a generating equation, for the proof see Proposition 2.2 in \cite{BarskiLochowski}. For their formulations we need a family of measures related to the pair $(G,Z)$. For $x\geq 0$ we define the measure
$$
\nu_{G(x)}(A):=\nu \{y \in \R^d: \langle G(x),y\rangle\in A \} , \quad A\in\mathcal{B}(\mathbb{R}).
$$
which is the image of the L\'evy measure $\nu$ under the linear transformation $y\mapsto \langle G(x), y\rangle$.
This measure may have an atom at zero and therefore its restriction $\nu_{G(x)}(\dd v)\mid_{v\neq 0}$  is used below.
The aforementioned conditions are as follows.

\begin{itemize}
\item The drift is affine 
\begin{gather}\label{linear drift}
F(x)=ax+b, \ \text{where} \ a\in\mathbb{R}, \ b\geq\int_{(1,+\infty)}(v-1)\nu_{G(0)}(\dd v).
\end{gather}
\item The covariance matrix of the Wiener part of $Z$ satisfies
\begin{align}\label{mult. CIR condition}
\frac{1}{2}\langle Q G(x), G(x)\rangle&=cx, \quad x\geq 0,
\end{align}
with some $c\geq 0$.

\item The jumps of $Z$ and the function $G$ are such that 
\begin{gather}\label{Z^G positive jumps}
\langle G(x), \triangle Z(t)\rangle\geq 0, \quad x \ge 0, t\geq 0,
\end{gather}
\begin{gather}\label{nu G0 finite variation}
\nu_{G(0)}(\dd v)=m(\dd v) \ \text{and} \ \int_{(0,+\infty)}v \ \nu_{G(0)}(\dd v)<+\infty,
\end{gather}
\begin{gather}\label{war calkowe na mu}
\int_{(0,+\infty)}(v \wedge v^2)\mu(\dd v)<+\infty,
\end{gather}
\begin{align}\label{rozklad nu G(x)}
\nu_{G(x)}(\dd v)\mid_{(0,+\infty)}&=\nu_{G(0)}(\dd v)\mid_{(0,+\infty)}+x\mu(\dd v), \quad x\geq 0.
\end{align}
\end{itemize}

Moreover, \eqref{generator Filipovica} reads
\begin{align}\label{generator R w tw}\nonumber
\mathcal{A}f(x)=cx f^{\prime\prime}(x)&+\Big[ax +b+\int_{(1,+\infty)}(1 -v)\{\nu_{G(0)}(\dd v)+x\mu(\dd v)\}\Big]f^{\prime}(x)\\[1ex]
&+\int_{(0,+\infty)}[f(x+v)-f(x)-f^{\prime}(x)(1\wedge v)]\{\nu_{G(0)}(\dd v)+x\mu(\dd v)\}.
\end{align}
In particular, with the parameters $a,b,c$ and the measures $\nu_{G(0)}(\dd v), \mu(\dd v)$ from \eqref{generator R w tw} at hand one can determine the zero coupon bond prices, for details see \cite{FilipovicATS}.

Note stronger integrability requirements  \eqref{nu G0 finite variation}, \eqref{war calkowe na mu} for the measures $m(\dd v), \mu(\dd v)$ than in \eqref{warunki na iary Filipovica}. They appear due to the fact that $Z$ is a martingale.

\begin{rem}\label{rem o rownaniu dla Laplasuf}
Conditions \eqref{mult. CIR condition}-\eqref{rozklad nu G(x)} describe the law of the family of one-dimensional L\'evy processes $Z^{G(x)}(t):=\langle G(x),Z(t)\rangle, x\geq 0$. Conditions \eqref{mult. CIR condition} and \eqref{rozklad nu G(x)} can be reformulated 
in terms of their Laplace exponents 
\begin{gather}\label{wqrunek z Laplacami}
J_{Z^{G(x)}}(b)=J_{Z}(b G(x))=cb^2+J_{\nu_{G(0)}}(b)+xJ_{\mu}(b), \quad b\geq 0,
\end{gather}
where $J_{\mu}(b):=\int_{0}^{+\infty}(e^{-bv}-1+bv)\mu(\dd v)$ and $J_{\nu_{G(0)}}$ is defined analogously.
\end{rem}

\begin{rem}\label{rem o charakterystykach dla alfa stab. 1 dim}
For the equation \eqref{rownanie na R stabilny d=1} with $\alpha\in(1,2)$ one can show that
$$
c=0, \ \nu_{G(0)}=0, \ \mu(\dd v)=\mathbf{1}_{\{v>0\}}\frac{1}{v^{1+\alpha}} \dd v,
$$
\end{rem}
hence $\mu(\dd v)$ is an $\alpha$-stable measure, for details see \cite{BarskiLochowski} or \cite{BarskiZabczyk}.

We start with an example of \eqref{main equation} where  $Z$ is an  $\alpha$-stable martingale in $\mathbb{R}^d, d>1$, with $\alpha\in(1,2)$. Recall, its radial measure is given by \eqref{radialna stabilna}. Since $Z$ has no Wiener part, the Laplace exponent of the jump part $X$ of $Z$ is identical with the Laplace exponent of $Z$ and admits the following representation:
\begin{align}\nonumber
J_X(z)&=\int_{\mathbb{S}^{d-1}}\lambda(\dd \xi)\int_{0}^{+\infty}\left(e^{-\langle z,r\xi\rangle}-1+\langle z,r\xi\rangle\right)\frac{1}{r^{1+\alpha}}\dd r\\[1ex]\nonumber
&=\int_{\mathbb{S}^{d-1}}\lambda(\dd \xi)\int_{0}^{+\infty}\left(e^{-r\langle z,\xi\rangle}-1+r\langle z,\xi\rangle\right)\frac{1}{r^{1+\alpha}}\dd r\\[1ex]\label{J dla alfa stabilnego}
&=c_{\alpha}\int_{\mathbb{S}^{d-1}}\langle z,\xi\rangle^\alpha\lambda(\dd \xi),
\end{align}
where $c_{\alpha}:=\Gamma(2-\alpha)/(\alpha(\alpha-1))$ and $\Gamma$ stands for the Gamma function.
Above we used the formula
$$
\int_{0}^{+\infty}\Big(e^{-uv}-1+uv\Big)\frac{1}{v^{1+\alpha}}\dd v=c_{\alpha} u^\alpha.
$$

\noindent
In the following example, assuming that $Z$ is an $\alpha$-stable martingale in $\mathbb{R}^{d}$,
we compute the condition for  the function $G$ in \eqref{main equation} so that  this equation 
generates an affine model.
\begin{ex}\label{ex1}
Let  $Z$ be an $\alpha$-stable martingale in $\mathbb{R}^d$ with the
Laplace exponent \eqref{J dla alfa stabilnego} and $G: [0,+\ns) \ra [0,+\ns)^d$, $G(0)=0$. Then the equation
\begin{gather}\label{rownanie z alfa stabilnym wiel.}
\dd R(t)=(aR(t)+b)\dd t+\langle G(R(t-)),Z(t)\rangle, 
\end{gather}
with $a\in\mathbb{R}, b\geq 0$ generates an affine model if and only if the function $G$ satisfies
\begin{gather}\label{war na G alfa stabilny}
\int_{\mathbb{S}^{d-1}}\langle G(x),\xi\rangle^\alpha \lambda(\dd \xi)=\frac{C}{c_{\alpha}}x,\quad x\geq 0,
\end{gather}
with $C\ge 0$. To prove this fact we need to show that \eqref{war na G alfa stabilny} is equivalent to 
\eqref{wqrunek z Laplacami} with some measure $\mu(\dd v)$. Since $Z$ has no Wiener part and 
$\nu_{G(0)}(\dd v)\equiv 0$, we see that \eqref{wqrunek z Laplacami} takes the form
$$
J_{Z}(bG(x))=J_{X}(bG(x))=xJ_{\mu}(b), \quad x,b\geq 0.
$$
By \eqref{J dla alfa stabilnego} 
$$
J_X(bG(x))=c_{\alpha}\int_{\mathbb{S}^{d-1}}\langle bG(x),\xi\rangle^\alpha\lambda(\dd \xi)
=c_{\alpha} b^{\alpha}\int_{\mathbb{S}^{d-1}}\langle G(x),\xi\rangle^\alpha\lambda(\dd \xi).
$$
Consequently,
$$
c_{\alpha} b^{\alpha}\int_{\mathbb{S}^{d-1}}\langle G(x),\xi\rangle^\alpha\lambda(\dd \xi)=xJ_{\mu}(b),
$$
holds if and only if
$$
J_\mu(b)=C b^{\alpha}, \quad \int_{\mathbb{S}^{d-1}}\langle G(x),\xi\rangle^\alpha\lambda(\dd \xi)=\frac{C}{c_{\alpha}}x,
$$
for some $C \ge 0$. Hence, $\mu$ is an $\alpha$-stable measure and $G$ can be any function satisfying \eqref{war na G alfa stabilny}. 
It follows from Remark \ref{rem o charakterystykach dla alfa stab. 1 dim} and \eqref{generator R w tw} that the generators of 
\eqref{rownanie z alfa stabilnym wiel.} and that of the equation \eqref{rownanie na R stabilny d=1}
are identical, so are the related bond markets. 
\end{ex}

\section{Reducibility of equations with multivariate noise}\label{sec Reducibility of equations with multivariate noise}

In this section we specify conditions for the equation \eqref{main equation}, written now  for convenience in the form
\begin{gather}\label{main equation 2}
\dd R(t)=(aR(t)+b)\dd t+\langle G(R(t-)), \dd Z(t)\rangle , \quad R(0)=x\geq 0,\quad t>0,
\end{gather}
to have the reducibility property.  The affine form of drift above is justified by \eqref{linear drift}. This means that \eqref{main equation 2} is supposed to generate
the same bond prices  as the equation
\begin{gather}\label{rownanie na R stabilny d=1 2}
\dd R(t)=(aR(t)+b)\dd t+C\cdot R(t-)^{1/\alpha} \dd Z^{\alpha}(t), 
\end{gather}
with $a\in\mathbb{R}, b\geq 0,C>0$ and an $\alpha$-stable real valued L\'evy process $Z^{\alpha}$ with some $\alpha\in(1,2)$. 
Recall that from Example \ref{ex1} we know that each generating equation \eqref{main equation 2} with $Z$ being an $\alpha$-stable process in $\mathbb{R}^d$ has the reducibility property. 

In \eqref{main equation 2}, $G:\mathbb{R}_{+}\longrightarrow \mathbb{R}^d$ and $Z$ is a L\'evy process and martingale in $\mathbb{R}^d$, called a {\it L\'evy martingale} for short.  It is characterized by a covariance matrix $Q$ of the Wiener part and a L\'evy measure $\nu$ which is assumed to admit a spherical decomposition
\begin{equation} \label{spherical 2}
\nu(A)=\int_{\S^{d-1}}\int_{0}^{+\infty}\mathbf{1}_{A}(r\xi)\gamma_{\xi}(\dd r) \ \lambda(\dd \xi), \quad A\in\mathcal{B}(\mathbb{R}^d),
\end{equation}
with a finite spherical measure $\lambda$  on the unit sphere $\S^{d-1}$ and some radial measures $\{\gamma_{\xi};\xi\in\S^{d-1}\}$. To avoid technical complications we can assume, and we do, the non-degeneracy condition for the radial measures, i.e.
\begin{gather}\label{non-degeneracy radial measures}
\xi\in \text{supp}(\lambda) \Longrightarrow \gamma_{\xi}\neq 0.
\end{gather}
If \eqref{non-degeneracy radial measures} is not satisfied, one can modify $\lambda$ by cutting off the part of its support 
where the radial measures disappear. This operation clearly does not affect \eqref{spherical 2}. Since $Z$ is a martingale, it follows from \eqref{Levy measure for a martingale} that
$$
\int_{\mathbb{R}^d}(\left| y\right|^2\wedge \left| y\right|)\nu(\dd y)=\int_{\S^{d-1}}\int_{0}^{+\infty}(\left|r\xi\right|^2\wedge \left|r\xi\right|)\gamma_\xi(\dd r)\lambda(\dd \xi)<+\infty,
$$
which means that 
\begin{gather}\label{levy measure condition martingale spherical}
\int_{0}^{+\infty}(r^2\wedge r)\gamma_\xi(\dd r)<+\infty, \quad \xi\in\text{supp}(\lambda).
\end{gather}
If the jump part of $Z$ has infinite variation, then it follows from \eqref{Levy measure finite variation} that 
\begin{gather}\label{infinite variation spherical}
\int_{\mid y\mid\leq 1}\left| y\right|\nu(dy)=\int_{\S^{d-1}} \int_{0}^{1}r \ \gamma_{\xi}(\dd r) \ \lambda(\dd \xi)=+\infty.
 \end{gather}
 We consider a stronger condition than \eqref{infinite variation spherical}, namely that
 \begin{gather}\label{warunek na miare radialna b}
\lambda(\Gamma_{\lambda})>0, \ \text{where} \ \Gamma_{\lambda}:=\left\{\xi\in \text{supp}(\lambda): \int_{0}^{1}r \gamma_{\xi}(dr)=+\infty\right\}.
\end{gather}
Consequently, if we assume that $\Gamma_{\lambda}$ is not contained in any proper linear subspace of $\R^d$, i.e.
\begin{gather}\label{infinite variation non-degenerated}
\text{Linear span}  \ (\Gamma_{\lambda}) = \R^d,
\end{gather}
then we obtain that
\begin{gather}\label{G0 znowu zero}
G(0)=0.
\end{gather}
To see this  let us notice that, by \eqref{Z^G positive jumps} and \eqref{spherical 2},
$
\lambda \cbr{\xi \in \S^{d-1}: \langle G(0), \xi \rangle <0} = 0
$
which implies that 
\begin{equation} \label{nonzerol}
\langle G(0),\xi \rangle \ge 0 \ \text{for any } \xi \in \text{supp} \ \lambda.
\end{equation}
By \eqref{nu G0 finite variation} we have 
\begin{align*}
\int_{0}^{+\infty}v \ \nu_{G(0)}(\dd v)&=\int_{\mathbb{R}^d}\langle G(0), y\rangle\nu(\dd y)\\[1ex]
&=\int_{\S^{d-1}}\langle G(0),\xi\rangle\int_{0}^{+\infty}r \ \gamma_{\xi}(\dd r) \lambda(\dd \xi)<+\infty,
\end{align*}
which, in view of \eqref{warunek na miare radialna b}, \eqref{infinite variation non-degenerated} and \eqref{nonzerol}
implies \eqref{G0 znowu zero}. Obviously, \eqref{infinite variation non-degenerated} also implies that
\begin{gather}\label{spherical part non-degenerated}
\text{Linear span  (supp } (\lambda)) = \R^d.
\end{gather}

\subsection{Main results}

For $\xi\in\text{supp}(\lambda)$ let us consider the Laplace exponent related to the measure $\gamma_{\xi}$, i.e.
$$
J_{\gamma_\xi}(b)=\int_{0}^{+\infty}(e^{-br}-1+br) \ \gamma_{\xi}(\dd r), \quad b\geq 0.
$$
We need the condition that there exists $K\geq 1$ such that 
\begin{gather}\label{warunek na miary radialne}
\sup_{\xi\in\text{supp}(\lambda)}J_{\gamma_\xi}(b)\leq K \cdot \inf_{\xi\in\text{supp}(\lambda)}J_{\gamma_\xi}(b), \quad b\geq 0.
\end{gather}

The main result of the paper is the following theorem.

\begin{tw} \label{tw_sferyczne}
Let $Z$ be a L\'evy martingale with a covariance matrix $Q$ of the Wiener part  and a L\'evy measure $\nu$ 
admitting the decomposition \eqref{spherical 2} with a spherical measure $\lambda$ satisfying \eqref{spherical part non-degenerated} and radial measures $\{\gamma_{\xi}; \xi\in{\S}^{d-1}\}$ satisfying \eqref{levy measure condition martingale spherical} and \eqref{warunek na miary radialne}. Let us also assume that \eqref{warunek na miare radialna b} and \eqref{infinite variation non-degenerated} are satisfied or that \eqref{G0 znowu zero} holds. Moreover, let $G:[0,+\infty)\longrightarrow \mathbb{R}^d$ be a continuous function such that  
\begin{gather}\label{istnienie G0}
G_0:=\lim_{x \ra 0^+} \frac{G(x)}{|G(x)|},
\end{gather} 
exists. 

Then if \eqref{main equation 2} generates an affine model, then the measure $\mu(\dd v)$ in \eqref{war calkowe na mu}-\eqref{rozklad nu G(x)} is $\alpha$-stable with $\alpha \in (1,2)$.
\end{tw}

The proof of Theorem \ref{tw_sferyczne} is presented in Subsection \ref{proof_of_tw_sf} and is preceded by 
some auxiliary results presented in Subsection \ref{Auxilliary results}.

\vskip 1ex

From Theorem \ref{tw_sferyczne} the following corollary and remarks follow.

\begin{cor} \label{cor_sferyczne} Let the assumptions of Theorem \ref{tw_sferyczne} be satisfied. If 
 \eqref{main equation 2} is a generating equation then \eqref{mult. CIR condition} is satisfied with $c=0$. This means that
 the continuous (Wiener) part of the process $Z^{G(x)}(t):=\langle G(x), Z(t)\rangle, t\geq 0$ vanishes for all $x >0$.
\end{cor}
{\bf Proof:}  It follows from \eqref{mult. CIR condition} that the Wiener part of $Z^{G(x)}$ satisfies
\begin{equation} \label{warrr4}
\frac{1}{2}\langle QG(x),G(x)\rangle=cx, \ x\geq 0 \ \text{for some} \  c\geq 0.
\end{equation}
Either directly by assumption \eqref{G0 znowu zero} or by the assumptions \eqref{warunek na miare radialna b} and \eqref{infinite variation non-degenerated} we get that $G(0)=0$. By Theorem \ref{tw_sferyczne}, the Laplace transform of the jump part of $Z$ satisfies 
\begin{equation} \label{warrr3}
J_{X}(b G(x))=J_{\nu_{G(0)}}(b)+xJ_{\mu}(b) = \gamma x b^{\alpha}, \quad x\geq 0\ \text{for some}\ \gamma > 0, \alpha \in (1,2).
\end{equation}
Condition \eqref{Z^G positive jumps} guarantees that $\langle G_0,y \rangle \ge 0$ for any $ y \in \text{supp } \nu$ and condition
\eqref{spherical part non-degenerated} guarantees that $y \mapsto \langle G_0, y\rangle, y\in  \text{supp} \ \nu$, does not vanish, hence $J_X(G_0) > 0$. Consequently, from \eqref{warrr3} we obtain
\[
\lim_{x \ra 0+} \frac{\gamma x}{|G(x)|^{\alpha}} =  \lim_{x \ra  0+} J_{X}\rbr{\frac{G(x)}{|G(x)|}} = J_{X}\rbr{G_{0}} \in (0, +\ns).
\]
From this, $\lim_{x \ra 0+} {|G(x)|} = 0$ and from \eqref{warrr4}  we further have 
\[
\langle Q G_{0},G_0\rangle = \lim_{x \ra 0+} \frac{\langle Q G(x),G(x)\rangle}{|G(x)|^2} =  \lim_{x \ra 0+} \frac{\gamma x}{|G(x)|^{\alpha}}  \frac{2c / \gamma}{|G(x)|^{2 - \alpha}} =
\begin{cases} 
0  \ \ \ \ \ \text{if } c = 0;\\
+ \ns   \ \text{if } c>0.
\end{cases} 
\]
Since  $\langle Q G_{0},G_0\rangle \neq +\ns$, we necessarily have  $c = 0$ which, in view of \eqref{warrr4}, means that the continuous (Wiener) part of $Z^{G(x)}$ vanishes.
\hfill $\square$

\begin{cor}\label{cor reducibility} 
It follows from Theorem \ref{tw_sferyczne} and Corollary \ref{cor_sferyczne} that each generating equation \eqref{main equation 2} satisfying assumptions of Theorem \ref{tw_sferyczne} satisfies conditions 
 \eqref{mult. CIR condition}-\eqref{rozklad nu G(x)} with
$$
c=0, \ \nu_{G(0)}=0, \ \mu(\dd v)=\mathbf{1}_{\{v>0\}}\frac{1}{v^{1+\alpha}} \dd v, \alpha\in(1,2).
$$
In view of Remark \ref{rem o charakterystykach dla alfa stab. 1 dim} the generators of \eqref{main equation 2} and of \eqref{rownanie na R stabilny d=1 2} are identical, so \eqref{main equation 2} has the reducibility property.
\end{cor}

\begin{rem}
In the formulation of Theorem \ref{tw_sferyczne} the assumption \eqref{istnienie G0} can be replaced by the existence 
of the limit $\lim_{x \ra +\ns} \frac{G(x)}{|G(x)|}$. Under the latter condition, however, we were unable to prove Corollary \ref{cor_sferyczne}.
\end{rem}

The following result provides sufficient conditions for the condition  \eqref{warunek na miary radialne} to hold.

\begin{prop}\label{prop o ograniczeniu miar promieniowych}
Let $\gamma(\dd r)$ and $\Gamma(\dd r)$ be two measures on $(0,+\infty)$ such that for any $\xi\in\text{supp}(\lambda)$
\begin{gather}\label{dominacja miar}
\gamma(A)\leq \gamma_\xi(A)\leq \Gamma(A), \quad A\in\mathcal{B}((0,+\infty)),
\end{gather}
and
\begin{gather}\label{warunki na miar ogr. z minimum}
0<\int_{0}^{+\infty}(r^2\wedge r)\ \gamma(\dd r)\leq\int_{0}^{+\infty}(r^2\wedge r)\ \Gamma(\dd r)<+\infty.
\end{gather}
If $\rbr{ \int_{\varepsilon}^1 r \gamma(\dd r) } \wedge \rbr{\int_1^{1/\varepsilon} r^2 \gamma(\dd r) }>0$ for all $\varepsilon >0$ sufficiently close to $0$ 
and there exist the limits 
\begin{gather*}
q_0:=\limsup_{\varepsilon \rightarrow 0+}\frac{\int_{\varepsilon}^{1}r\Gamma(\dd r)}{\int_{\varepsilon}^{1}r\gamma(\dd r)},\qquad 
q_{\infty}:=\limsup_{\varepsilon \rightarrow 0+}\frac{\int_{1}^{1/\varepsilon}r^2\Gamma(\dd r)}{\int_{1}^{1/\varepsilon}r^2\gamma(\dd r)},
\end{gather*}
and both are finite, then \eqref{warunek na miary radialne} is satisfied.
\end{prop}

\begin{ex}[Spherically balanced L\'evy measure] Let us consider 
the case when the radial measures satisfy
\begin{gather*}
\gamma(A)\le\gamma_{\xi}(A)\le K\cdot\gamma(A), \quad A\in\mathcal{B}((0,+\infty)),
\end{gather*}
with some finite constant $K\geq 1$ and a measure $\gamma$ such that
$$
0<\int_{0}^{+\infty}(r^2\wedge r)\ \gamma(\dd r)<+\infty
$$
and $\rbr{ \int_{\varepsilon}^1 r \gamma(\dd r) } \wedge \rbr{\int_1^{1/\varepsilon} r^2 \gamma(\dd r) }>0$ for all $\varepsilon >0$ sufficiently close to $0$.
Then $q_0\le K$ and $q_{\infty} \le K$, so by Proposition \ref{prop o ograniczeniu miar promieniowych} condition \eqref{warunek na miary radialne} is satisfied.
\end{ex}

\noindent
{\bf Proof of Proposition \ref{prop o ograniczeniu miar promieniowych}:} Under \eqref{dominacja miar} we clearly have 
$$
J_{\gamma}(b)\leq \inf_{\xi\in\text{supp}(\lambda)}J_{\gamma_{\xi}}(b)  \leq\sup_{\xi\in\text{supp}(\lambda)}J_{\gamma_{\xi}}(b)\leq J_{\Gamma}(b), \quad b\geq 0,
$$
where
\begin{gather*}
J_{\gamma}(b):=\int_{0}^{+\infty}(e^{-br}-1+br) \ \gamma(\dd r), \quad J_{\Gamma}(b):=\int_{0}^{+\infty}(e^{-br}-1+br) \ \Gamma(\dd r), \quad b\geq 0.
\end{gather*}
Therefore \eqref{warunek na miary radialne} is implied by the condition
\begin{gather}\label{warunek zamienny}
J_{\Gamma}(b)\leq K \cdot J_{\gamma}(b), \quad b\geq 0.
\end{gather}
Since the functions $J_{\gamma}(\cdot), J_{\Gamma}(\cdot)$ are continuous, hence bounded on compacts,  \eqref{warunek zamienny} is satisfied with some $K\geq 1$ if and only if
\begin{equation}
p_{\infty}:=\limsup_{b\rightarrow+\infty}\frac{J_{\Gamma}(b)}{J_{\gamma}(b)}<+\infty \label{eq:mmain}
\end{equation}
and
\begin{equation}
p_{0}:=\limsup_{b\rightarrow0+}\frac{J_{\Gamma}(b)}{J_{\gamma}(b)}<+\infty. \label{eq:mmmain}
\end{equation}
In what follows we show that \eqref{eq:mmain} and  \eqref{eq:mmmain} indeed hold.

Let us notice that for $x\ge0$
\[
e^{-x}-1+x\sim x\wedge x^{2},
\]
where the relation $\sim$ means that there exist universal positive
numbers $k$ and $K$ such that
\[
k\cdot x\wedge x^{2}\le e^{-x}-1+x\le K\cdot x\wedge x^{2}.
\]
Thus, to prove (\ref{eq:mmain}) it is sufficient to prove that
\begin{align*}
 & \limsup_{b\rightarrow+\infty}\frac{\int_{0}^{+\infty}\min(br,b^{2}r^{2})\Gamma(\mathrm{d}r)}{\int_{0}^{+\infty}\min(br,b^{2}r^{2})\gamma(\mathrm{d}r)}\\
 & =\limsup_{b\rightarrow+\infty}\frac{\int_{0}^{1/b}b^{2}r^{2}\Gamma(\mathrm{d}x)+\int_{1/b}^{+\infty}br\Gamma(\mathrm{d}r)}{\int_{0}^{1/b}b^{2}r^{2}\gamma(\mathrm{d}x)+\int_{1/b}^{+\infty}br\gamma(\mathrm{d}r)}<+\infty.
\end{align*}
Let us define the functions
\[
G(y):=\int_{(y,+\infty)}r\Gamma(\mathrm{d}r),\quad g(y):=\int_{(y,+\infty)}r\gamma(\mathrm{d}r),\quad y>0.
\]
By integration by parts, 
\begin{equation}
\int_{0}^{1/b}r^{2}\Gamma(\mathrm{d}r)=\int_{0}^{1/b}r(-\mathrm{d}G(r))=-r\cdot G(r)|_{0}^{1/b}+\int_{0}^{1/b}G(r)\mathrm{d}r.\label{eq:one_step}
\end{equation}
We fix $\varepsilon>0$ and for $y\in(0,\varepsilon)$ estimate 
\begin{align*}
y\cdot G(y) & =y\int_{y}^{+\infty}r\Gamma(\mathrm{d}r) = \int_{y}^{\varepsilon}yr\Gamma(\mathrm{d}r)+y\cdot G(\varepsilon)\\
 & \le\int_{0}^{\varepsilon}r^{2}\Gamma(\mathrm{d}r)+y\cdot G(\varepsilon).
\end{align*}
From this it follows 
\[
\limsup_{y\rightarrow0+}y\cdot G(y)\le\int_{0}^{\varepsilon}r^{2}\Gamma(\mathrm{d}r)
\]
and by the finiteness of $\int_{0}^{1}r^{2}\Gamma(\mathrm{d}r)$ and
arbitrary choice of $\varepsilon$, we get 
\[
\limsup_{y\rightarrow0+}yG(y)=0.
\]
Thus, (\ref{eq:one_step}) takes the form
\[
\int_{0}^{1/b}r^{2}\Gamma(\mathrm{d}r)=\int_{0}^{1/b}r(-\mathrm{d}G(r))=-\frac{1}{b}\cdot G\left(\frac{1}{b}\right)+\int_{0}^{1/b}G(r)\mathrm{d}r.
\]
Similarly, 
\[
\int_{0}^{1/b}r^{2}\gamma(\mathrm{d}r)=-\frac{1}{b}\cdot g\left(\frac{1}{b}\right)+\int_{0}^{1/b}g(r)\mathrm{d}r.
\]
Now we calculate
\begin{align*}
 & \limsup_{b\rightarrow+\infty}\frac{\int_{0}^{1/b}b^{2}r^{2}\Gamma(\mathrm{d}x)+\int_{1/b}^{+\infty}br\Gamma(\mathrm{d}r)}{\int_{0}^{1/b}b^{2}r^{2}\gamma(\mathrm{d}x)+\int_{1/b}^{+\infty}br\gamma(\mathrm{d}r)}\\
 & =\limsup_{b\rightarrow+\infty}\frac{b^{2}\left(-\frac{1}{b}\cdot G\left(\frac{1}{b}\right)+\int_{0}^{1/b}G(r)\mathrm{d}r\right)+b\cdot G\left(\frac{1}{b}\right)}{b^{2}\left(-\frac{1}{b}\cdot g\left(\frac{1}{b}\right)+\int_{0}^{1/b}g(r)\mathrm{d}r\right)+b\cdot g\left(\frac{1}{b}\right)}\\
 & =\limsup_{b\rightarrow+\infty}\frac{\int_{0}^{1/b}G(r)\mathrm{d}r}{\int_{0}^{1/b}g(r)\mathrm{d}r}<+\infty,
\end{align*}
where the last estimate follows from the assumption
\[
q_0 = \limsup_{y\rightarrow0+}\frac{\int_{y}^{1}r\Gamma(\mathrm{d}r)}{\int_{y}^{1}r\Gamma(\mathrm{d}r)}<+\infty
\]
and the finiteness of $\int_{1}^{+\infty}r\Gamma(\mathrm{d}r),$ which
yields that the ratio $G(r)/g(r)$ is separated from $+\infty$ for
$r$ sufficiently close to $0$.

To prove (\ref{eq:mmmain}) we will proceed in a similar way. To prove
(\ref{eq:mmain}) it is sufficient to prove that
\begin{align*}
 & \limsup_{b\rightarrow0+}\frac{\int_{0}^{+\infty}\min(br,b^{2}r^{2})\Gamma(\mathrm{d}r)}{\int_{0}^{+\infty}\min(br,b^{2}r^{2})\gamma(\mathrm{d}r)}\\
 & =\limsup_{b\rightarrow0+}\frac{\int_{0}^{1/b}b^{2}r^{2}\Gamma(\mathrm{d}x)+\int_{1/b}^{+\infty}br\Gamma(\mathrm{d}r)}{\int_{0}^{1/b}b^{2}r^{2}\gamma(\mathrm{d}x)+\int_{1/b}^{+\infty}br\gamma(\mathrm{d}r)}<+\infty.
\end{align*}
We define 
\[
Q(y):=\int_{(0,y]}r^{2}\Gamma(\mathrm{d}r),\quad q(y):=\int_{(0, y]}r^{2}\gamma(\mathrm{d}r),\quad y>0.
\]
By integration by parts, 
\begin{equation}
\int_{1/b}^{+\infty}r\Gamma(\mathrm{d}r)=\int_{1/b}^{+\infty}\frac{1}{r}\mathrm{d}Q(r)=\frac{1}{r}Q(r)|_{1/b}^{+\infty}+\int_{1/b}^{+\infty}\frac{Q(r)}{r^{2}}\mathrm{d}r.\label{eq:one_step-1}
\end{equation}
We fix $M>0$ and for $y\in(M,+\infty)$ estimate 
\begin{align*}
\frac{1}{y}Q(y) & =\frac{1}{y}\int_{0}^{y}r^{2}\Gamma(\mathrm{d}r) = \frac{1}{y}\int_{0}^{M}r^{2}\Gamma(\mathrm{d}r)+\int_{M}^{y}\frac{r}{y}r\Gamma(\mathrm{d}r)\\
 & \le\frac{1}{y}Q(M)+\int_{M}^{+\ns}r\Gamma(\mathrm{d}r).
\end{align*}
From this it follows 
\[
\limsup_{y\rightarrow+\infty}\frac{1}{y}Q(y)\le\int_{M}^{+\ns}r\Gamma(\mathrm{d}r)
\]
and by the finiteness of $\int_{1}^{+\infty}r\Gamma(\mathrm{d}r)$
and arbitrary choice of $M$, we get 
\[
\limsup_{y\rightarrow+\infty}\frac{1}{y}Q(y)=0.
\]
Thus, (\ref{eq:one_step-1}) takes the form
\[
\int_{1/b}^{+\infty}r\Gamma(\mathrm{d}r)=-bQ\left(\frac{1}{b}\right)+\int_{1/b}^{+\infty}\frac{Q(r)}{r^{2}}\mathrm{d}r.
\]
Similarly, 
\[
\int_{0}^{1/b}r^{2}\gamma(\mathrm{d}r)=-bq\left(\frac{1}{b}\right)+\int_{1/b}^{+\infty}\frac{q(r)}{r^{2}}\mathrm{d}r.
\]
Now we calculate
\begin{align*}
 & \limsup_{b\rightarrow0+}\frac{\int_{0}^{1/b}b^{2}r^{2}\Gamma(\mathrm{d}x)+\int_{1/b}^{+\infty}br\Gamma(\mathrm{d}r)}{\int_{0}^{1/b}b^{2}r^{2}\gamma(\mathrm{d}x)+\int_{1/b}^{+\infty}br\gamma(\mathrm{d}r)}\\
 & =\limsup_{b\rightarrow0+}\frac{b^{2}Q\left(\frac{1}{b}\right)+b\left(-bQ\left(\frac{1}{b}\right)+\int_{1/b}^{+\infty}Q(r)\frac{\mathrm{d}r}{r^{2}}\right)}{b^{2}q\left(\frac{1}{b}\right)+b\left(-bq\left(\frac{1}{b}\right)+\int_{1/b}^{+\infty}q(r)\frac{\mathrm{d}r}{r^{2}}\right)}\\
 & =\limsup_{b\rightarrow0+}\frac{\int_{1/b}^{+\infty}Q(r)\frac{\mathrm{d}r}{r^{2}}}{\int_{1/b}^{+\infty}q(r)\frac{\mathrm{d}r}{r^{2}}}<+\infty,
\end{align*}
where the last estimate follows from the assumption
\[
q_{\ns} = \limsup_{y\rightarrow0+}\frac{\int_{1}^{1/y}r^{2}\Gamma(\mathrm{d}r)}{\int_{1}^{1/y}r^{2}\Gamma(\mathrm{d}r)}<+\infty
\]
and the finiteness of $\int_{0}^{1}r^{2}\Gamma(\mathrm{d}r),$ which
yields that the ratio $Q(r)/q(r)$ is separated from $+\infty$ for
sufficiently large $r$.
 \hfill$\square$

\subsubsection{Jump measures with densities}\label{sec o Laplacach dla gestosci wiki}

In this subsection we formulate conditions required for the reducibility of \eqref{main equation 2} 
in the important case when the L\'evy measure of $Z$ has a density, i.e.
$\nu(\dd x)=g(x) \dd x$. Let us consider the polar transformation $\xi:P:=[0,\pi]^{d-2}\times[0,2\pi]\rightarrow {\mathbb{S}^{d-1}}$ given by 
\begin{align*}
\xi_1&=\cos\alpha_1,\\
\xi_2&=\sin\alpha_1\cdot\cos\alpha_2,\\
\xi_3&=\sin\alpha_1\cdot \sin\alpha_2\cdot\cos\alpha_3,\\
\vdots\\
\xi_{d-1}&=\sin\alpha_1\cdot \sin\alpha_2\cdot...\cdot\sin\alpha_{d-2}\cdot\cos\alpha_{d-1},\\
\xi_d&=\sin\alpha_1\cdot \sin\alpha_2\cdot...\cdot\sin\alpha_{d-2}\cdot\sin\alpha_{d-1}.
\end{align*}
The change of variables for polar coordinates $x=r\xi$ yields
\begin{align}\label{wzorek wiki}\nonumber
&\int_{\mathbb{R}^{d}}f(x)\nu(\dd x)=\int_{\mathbb{R}^{d}}f(x) g(x)\dd x\\
&=\int_{P} \int_{0}^{+\infty}\left(f(r\xi)g(r\xi)\cdot r^{d-1}\sin^{d-2}\alpha_1\cdot\sin^{d-3}\alpha_2\cdot\sin^{d-4}\alpha_3\cdot...\cdot\sin\alpha_{d-2}\right) \ \dd r \ \dd \alpha_1...\dd\alpha_{d-1},
\end{align}
for a $\nu$-integrable function $f$. Noting that 
$$
\sin^{d-2}\alpha_1\cdot\sin^{d-3}\alpha_2\cdot\sin^{d-4}\alpha_3\cdot...\cdot\sin\alpha_{d-2} =
\sqrt{1-\xi_1^2}\cdot \sqrt{1-(\xi_1^2+\xi_2^2)}\cdot...\cdot
\sqrt{1-(\xi_1^2+...+\xi_{d-2}^2)}
$$
we write \eqref{wzorek wiki} in the form
\begin{gather}\label{rozkład sferyczny z gestoscia}
\int_{\mathbb{S}^{d-1}}\int_{0}^{+\infty}f(r\xi)g(r\xi) r^{d-1}\sqrt{1-\xi_1^2}\cdot \sqrt{1-(\xi_1^2+\xi_2^2)}\cdot...\cdot
\sqrt{1-(\xi_1^2+...+\xi_{d-2}^2)} \ \dd r \ \lambda(\dd \xi).
\end{gather}
Above $\lambda$ stands for the image of the Lebesgue measure on $P$ under the transformation $\xi:P\rightarrow \S^{d-1}$ restricted to the set 
\begin{gather}\label{mathcal G}
\mathcal{G}:=\{\xi\in\S^{d-1}: g(r\xi) \not\equiv 0, r\geq 0\}.
\end{gather}
This definition of $\lambda$ is in accordance with \eqref{non-degeneracy radial measures} and implies that 
$$
\text{supp}(\lambda)=\mathcal{G}.
$$
From \eqref{rozkład sferyczny z gestoscia} we see that the radial measures have the form
\begin{gather}\label{wzor na gamma wiki}
\gamma_{\xi}(\dd r)=g(r\xi) r^{d-1}\sqrt{1-\xi_1^2}\cdot \sqrt{1-(\xi_1^2+\xi_2^2)}\cdot...\cdot
\sqrt{1-(\xi_1^2+...+\xi_{d-2}^2)} \ \dd r,
\end{gather}
and we can reformulate conditions needed in Theorem  \ref{tw_sferyczne} in terms of $g$. It is clear that
in \eqref{levy measure condition martingale spherical} we can use the equivalence, for $\xi\in\text{supp}(\lambda)$,
$$
 \int_{0}^{1}r\gamma_{\xi}(\dd r)=+\infty \quad \Longleftrightarrow \quad  \int_{0}^{1}r^{d}g(r\xi)\dd r=+\infty.
$$

From Theorem \ref{tw_sferyczne}, Corollary \ref{cor_sferyczne}, Corrolary \ref{cor reducibility} and Proposition \ref{prop o ograniczeniu miar promieniowych} one can deduce the following result.
 
\begin{tw}\label{th reducibility density}
Let $Z$ be a L\'evy martingale with a covariance matrix $Q$ of the Wiener part  and a L\'evy measure with density 
$\nu(\dd x)=g(x)\dd x$ satisfying the following conditions
\begin{enumerate}[a)]
\item $\int_{\mathbb{R}^d}(\left| x\right|^2\wedge \left| x\right|) \ g(x)\dd x<+\infty$,
\item $\emph{Linear span}(\mathcal{G})=\mathbb{R}^d$ with $\mathcal{G}$ given by \eqref{mathcal G},
\item 
$\label{wiecej niz niesk. wahanie}
\lambda\left(\xi\in\mathcal{G}: \int_{0}^{1}r^{d}g(r\xi)\dd r=+\infty\right)>0.
$
\end{enumerate}
Let us define the functions
\begin{align}\label{sup inf definitions d>2 1}
\underline{g}(r)&:= \inf_{\left|x\right| = r}g(x)\sqrt{1-\frac{x_1^2}{\left|x\right|^2}}\cdot \sqrt{1-\frac{x_1^2+x_2^2}{\left|x\right|^2}}\cdot...\cdot
\sqrt{1-\frac{x_1^2+x_2^2+...+x_{d-2}^2}{\left|x\right|^2}}, \quad r\geq 0,\\\label{sup inf definitions d>2 2}
\bar{g}(r)&:= \sup_{\left|x\right| = r}g(x)\sqrt{1-\frac{x_1^2}{\left|x\right|^2}}\cdot \sqrt{1-\frac{x_1^2+x_2^2}{\left|x\right|^2}}\cdot...\cdot
\sqrt{1-\frac{x_1^2+x_2^2+...+x_{d-2}^2}{\left|x\right|^2}}, \quad r\geq 0,
\end{align}
and assume that  they satisfy
\begin{gather}\label{war cccalk d,d+1}
0<\int_{0}^{+\infty}(r^d\wedge r^{d+1}) \ \underline{g}(r)\dd r\leq \int_{0}^{+\infty}(r^d\wedge r^{d+1}) \ \bar{g}(r)\dd r<+\infty,
\end{gather}
and
\begin{gather}\label{warunki epsilonowe d>2}
\limsup_{\varepsilon\rightarrow 0+}\frac{\int_{\varepsilon}^{1}r^d \bar{g}(r)\dd r}{\int_{\varepsilon}^{1}r^d \underline{g}(r)\dd r}<+\infty \quad \text{and} \quad \limsup_{\varepsilon\rightarrow 0+}\frac{\int_{1}^{1/\varepsilon}r^{d+1} \bar{g}(r)\dd r}{\int_{1}^{1/\varepsilon}r^{d+1} \underline{g}(r)\dd r}<+\infty,
\end{gather}
and denominators in \eqref{warunki epsilonowe d>2} are positive for all $\varepsilon>0$ sufficiently close to $0$.
Let us also assume that $G:[0,+\infty)\longrightarrow \mathbb{R}^d$ is a continuous function such that  
\begin{gather*}\label{istnienie G0}
G_0:=\lim_{x \ra 0^+} \frac{G(x)}{|G(x)|},
\end{gather*} 
exists. 

Then, if \eqref{main equation 2} generates an affine model, then it has the reducibility property.
\end{tw}
 
In the proof of Theorem \ref{th reducibility density} one shows that \eqref{dominacja miar} is satisfied. Indeed, 
in view of \eqref{wzor na gamma wiki}, \eqref{sup inf definitions d>2 1} and \eqref{sup inf definitions d>2 2} we have
$$
\underline{g}(r) r^{d-1} \dd r \leq \gamma_{\xi}(\dd r)\leq \bar{g}(r) r^{d-1} \dd r.
$$
For the measures $\gamma(\dd r):=\underline{g}(r) r^{d-1} \dd r$ and $\Gamma(\dd r):=\bar{g}(r) r^{d-1} \dd r$ we recognize
in \eqref{war cccalk d,d+1}, \eqref{warunki epsilonowe d>2} the conditions \eqref{warunki na miar ogr. z minimum}, $q_0<+\infty$ and $q_{\infty}<+\infty$. It follows from Proposition \ref{prop o ograniczeniu miar promieniowych}  that condition  \eqref{warunek na miary radialne} is satisfied.

\begin{rem}
 In the case $d=2$ the functions $\underline{g}(r)$, $\bar{g}(r)$ take a simple form, i.e.
\begin{gather}\label{sup inf definitions}
\underline{g}(r) := \inf_{\left|x\right| = r}g(x), \qquad \bar{g}(r) := \sup_{\left|x\right| = r}g(x), 
\end{gather}
and therefore \eqref{war cccalk d,d+1} and \eqref{warunki epsilonowe d>2} in
Theorem \ref{th reducibility density} provide simple conditions for the reducibility of \eqref{main equation 2}.
\end{rem}

\subsection{Auxilliary results}\label{Auxilliary results}

Let us consider a generating equation \eqref{main equation 2} with some function $G$ and a L\'evy martingale $Z$ with the Laplace exponent of its jump part $J_X(\cdot)$, see \eqref{Jdef} for definition. Our first aim is to estimate the function $J_X(b G(x)), b,x\geq 0$ with the use of the function $J_X(b G_{0}), b\geq 0$ for $x$ such that $G(x)/\left|G(x)\right|$ is close to $G_0$. 
The solution of this problem is presented in Lemma \ref{lem oszacownaia na funkcje podcalkowa}, Proposition \ref{positivity} and Proposition \ref{sferycznosc}.

Let $\rho(\dd v)$ be an auxiliary L\'evy measure on $(0,+\infty)$ satisfying
\begin{gather}\label{Ro}
\int_{0}^{+\infty}(v^2\wedge v)\rho(\dd v)<+\infty,
\end{gather}
and 
\begin{gather}\label{funkcja Jot}
J_\rho(z):=\int_{(0,+\infty)}(e^{-zv}-1+zv)\rho(\dd v), \quad z\geq 0,
\end{gather}
The second aim of this section is to provide sufficient conditions for $J_\rho$
to be a power function. This problem is solved in Lemma \ref{lem o charakteryzacji f potegowej} and Lemma \ref{lem Weyl}.

\begin{lem}\label{lem oszacownaia na funkcje podcalkowa}
The function $H:[0,+\infty)\longrightarrow \mathbb{R}$ given by
$$
H(z)=e^{-z}-1+z,
$$
is convex, strictly increasing and 
\begin{gather}\label{min max}
\min\{1,t^2\}\cdot H(z)\leq H(t z)\leq\max\{1,t^2\}\cdot H(z),\quad z\geq 0,t>0.
\end{gather}
\end{lem}
{\bf Proof:} Since $H^\prime(z)=1-e^{-z}$ the monotonicity and convexity of $H$ follows. For $t\geq 1$ it follows from the monotonicity of $H$ that
$$
H(tz)\geq H(z)=\min\{1,t^2\}H(z).
$$
Let us notice that the function
$$
t \longrightarrow\frac{(1-e^{-t})t}{e^{-t}-1+t}, \quad t\geq 0,
$$
is strictly decreasing, with limit $2$ at zero and $1$ at infinity. This implies that
 \begin{equation} \label{oszH}
 (e^{-t}-1+t) < (1-e^{-t})t < 2 (e^{-t}-1+t), \quad t \in (0, +\ns).
 \end{equation}
From \eqref{oszH} we obtain
$$
\frac{\dd}{\dd s}\ln H(s)=\frac{H^\prime(s)}{H(s)}=\frac{1-e^{-s}}{e^{-s}-1+s}\leq \frac{2}{s}, \quad s>0,
$$
and, consequently, we obtain that for $t\geq 1$: 
$$
\ln H(tz)-\ln H(z)\leq \int_{z}^{tz}\frac{2}{s}\dd s=\ln t^2.
$$
Thus 
\begin{gather}\label{ktktk}
\min\{1,t^2\}H(z) = H(z)  \le H(tz)\leq t^2 H(z)=\max\{1,t^2\}H(z).
\end{gather}
Using the monotonicity of $H$ and  \eqref{ktktk} we see that for $t\in(0,1)$:
$$
H(tz)\leq H(z)=H\rbr{\frac{1}{t}tz}\leq \frac{1}{t^2}H(tz),
$$
so also for $t\in(0,1)$
$$
\min\{1,t^2\}H(z)=t^2 H(z)\leq H(tz)\leq H(z)=\max\{1,t^2\} H(z).
$$
\hfill $\square$

\begin{cor} \label{lem oszacownaia na wzrost transformaty}
It follows from \eqref{min max} and the formula
$$
J_\rho(z):=\int_{(0,+\infty)}H(zv)\rho(\dd v)<+\infty
$$
that the function $J_\rho$ satisfies
\begin{gather}\label{min maxxx}
\min\cbr{1,t^2}\cdot J_\rho(z)\leq J_\rho(t z)\leq\max\cbr{1,t^2}\cdot J_\rho(z),\quad z\geq 0,t>0.
\end{gather}
\end{cor}

\begin{prop} \label{positivity}
If \eqref{main equation 2} generates an affine model and $G_{\ns}$ is an arbitraty limit point of the set
\[
\cbr{\frac{G(x)}{\left|G(x)\right|}:x>0}
\]
then
\[
\nu\cbr{y\in\R^{d}:\left\langle G_{\ns},y\right\rangle <0}=0.
\]
\end{prop}
{\bf Proof:}  Assume that
\[
\nu\cbr{y\in\R^{d}:\left\langle G_{\ns},y\right\rangle <0}=\nu\cbr{y\in\R^{d}\setminus \{0\}:\left\langle G_{\ns},\frac{y}{\left|y\right|}\right\rangle <0}>0.
\]
Then there exists a natural $n$ such that for
\[
V_{n}:=\cbr{y\in\R^{d}\setminus\{0\}:\left\langle G_{\ns},\frac{y}{\left|y\right|}\right\rangle <-\frac{1}{n}}
\]
 one has $\nu\rbr{V_{n}}>0.$

Let $x$ be such that 
\[
\left|\frac{G(x)}{\left|G(x)\right|}-G_{\ns}\right|\le\frac{1}{2n}.
\]
It follows from the Schwarz inequality that for any $y\in\R^{d}$,
\begin{equation}
\left|\left\langle \frac{G(x)}{\left|G(x)\right|},y\right\rangle -\left\langle G_{\ns},y\right\rangle \right|\le\left|\frac{G(x)}{\left|G(x)\right|}-G_{\ns}\right|\left|y\right|\le\frac{1}{2n}\left|y\right|.\label{eq:schwarz}
\end{equation}
Let $y\in V_{n}$. From (\ref{eq:schwarz}) and the definition of
$V_{n}$ we estimate
\[
\left\langle \frac{G(x)}{\left|G(x)\right|},y\right\rangle \le\left\langle G_{\ns},y\right\rangle +\frac{1}{2n}\left|y\right|<-\frac{1}{n}\left|y\right|+\frac{1}{2n}\left|y\right|=-\frac{1}{2n}\left|y\right|<0.
\]
Hence 
\[
\nu\cbr{y\in\R^{d}:\left\langle \frac{G(x)}{\left|G(x)\right|},y\right\rangle <0}\ge\nu\rbr{V_{n}}>0
\]
which is a contradiction with \eqref{Z^G positive jumps}.
\hfill$\square$

\begin{prop} \label{sferycznosc} 
Let us assume that \eqref{main equation 2} is generating equation and that $\nu$ has the form \eqref{spherical 2} where $\lambda$ satisfies \eqref{spherical part non-degenerated} and $\gamma_\xi(\dd r)$ satisfies \eqref{warunek na miary radialne}. Let $G_{\ns}$ be any limit point of the set 
\[
\cbr{\frac{G(x)}{\left|G(x)\right|}:x>0}.
\]
Define
\[
M_{G_{\ns}}(b):=J_{X}\rbr{b\cdot G_{\ns}} =\int_{\S^{d-1}}\int_{0}^{+\ns}H\rbr{b\left\langle G_{\ns},r\cdot \xi\right\rangle }\gamma_{\xi}\rbr{\dd r}\lambda(\dd \xi),
\]
where $H(z):= e^{-z}-1+z$.
There exists a function $\delta: (0,1) \ra (0,+\ns)$ such that for any $\varepsilon_{0}>0$, any $b\ge0$ and $x>0$ such that $\left|\frac{G(x)}{\left|G(x)\right|}-G_{\ns}\right|\le \delta\rbr{\varepsilon_0}$ we have
 
\begin{equation}
\left(1-\varepsilon_{0}\right)M_{G_{\ns}}(b\left|G(x)\right|)\le J_{X}\rbr{bG(x)}\le\left(1+\varepsilon_{0}\right)M_{G_{\ns}}(b\left|G(x)\right|).\label{eq:teeza}
\end{equation}
\end{prop}
{\bf Proof:}  Let $\varepsilon\in(0,1)$ be such that 
\begin{equation}
\left(1+\varepsilon\right)^{2}\left(1+\frac{4 K \varepsilon}{\left(1-\varepsilon\right)^{3}}\right)\le1+\varepsilon_{0},\quad\frac{\left(1-\varepsilon\right)^{2}}{\left(1+\frac{ K \varepsilon}{1-\varepsilon}\right)}\ge1-\varepsilon_{0}.\label{eq:mniam_mniam_mniam}
\end{equation}
Let us assume that 
\begin{equation} \label{normalisation}
\lambda\cbr{\xi\in\S^{d-1}:\left\langle G_{\ns},\xi\right\rangle >0} = \lambda\rbr{\S^{d-1}}-\lambda\cbr{\xi\in\S^{d-1}:\left\langle G_{\ns},\xi\right\rangle =0}=1,
\end{equation}
(we can assume this, multiplying $\lambda$ by a positive constant, provided that \[\lambda\cbr{\xi\in\S^{d-1}:\left\langle G_{\ns},\xi\right\rangle >0}>0,\]
otherwise it follows from Proposition \ref{positivity} that we get a
degenerated case \[\lambda\rbr{\S^{d-1}}=\lambda\cbr{\xi\in\S^{d-1}:\left\langle G_{\ns},\xi\right\rangle =0}\] 
where \eqref{spherical part non-degenerated} is broken).
Let $\eta\in(0,1)$ be such that
\begin{equation} \label{normalisation1}
\lambda\cbr{\xi\in\S^{d-1}:0<\left\langle G_{\ns},\xi\right\rangle <\eta}\le\varepsilon.
\end{equation}
Moreover, by Proposition \ref{positivity}, 
\begin{align*}
0=\nu\cbr{y\in\R^{d}:\left\langle G_{\ns},y\right\rangle <0}&=
\int_{\S^{d-1}}\int_{0}^{+\infty}r\langle G_{\infty},\xi\rangle\gamma_{\xi}(\dd r)\lambda(\dd \xi)\\
&\geq\lambda\cbr{\xi\in\S^{d-1}:\left\langle G_{\ns},\xi\right\rangle <0}\cdot \sup_{\xi\in\S^{d-1}}\gamma_{\xi}\rbr{\R_{+}},
\end{align*}
so it follows that 
$$
\lambda\cbr{\xi\in\S^{d-1}:\left\langle G_{\ns},\xi\right\rangle <0}=0.
$$
Let us define 
\[
\mathbb{V}_{\eta}=\cbr{\xi\in\S^{d-1}:0<\left\langle G_{\ns},\xi\right\rangle <\eta}.
\]
Let $x$ be such that 
\[
\left|\frac{G(x)}{\left|G(x)\right|}-G_{\ns}\right|\le \delta\rbr{\varepsilon_0} :=  \eta\cdot\varepsilon.
\]
From Lemma \ref{lem oszacownaia na funkcje podcalkowa}, for $b,r\ge0$ and $\xi\in\S^{d-1}$ such that $\left\langle G_{\ns},\xi\right\rangle \in[0,\eta)$
we estimate 
\begin{align}\label{pierewasze szacownaiw}\nonumber
  H&\rbr{b\cdot r\left\langle G(x),\xi\right\rangle } \le H\rbr{b\cdot r\left|G(x)\right|\left(\left\langle G_{\ns},\xi\right\rangle +\left|\left\langle \frac{G(x)}{\left|G(x)\right|}-G_{\ns},\xi\right\rangle \right|\right)}\\\nonumber
 & \le\max\left\{H\rbr{b\cdot r\left|G(x)\right|2\left\langle G_{\ns},\xi\right\rangle },H\rbr{b\cdot r\cdot\left|G(x)\right|2\left|\left\langle \frac{G(x)}{\left|G(x)\right|}-G_{\ns},\xi\right\rangle \right|}\right\}\\\nonumber
 & \le\max\left\{H\rbr{b\cdot r\left|G(x)\right|2\eta},H\rbr{b\cdot r\cdot\left|G(x)\right|2\eta\cdot\varepsilon}\right\}\\\nonumber
 & =H\rbr{b\cdot r\left|G(x)\right|2\eta}\\
 & \le4H\rbr{b\cdot r\left|G(x)\right|\eta}.
\end{align}
It follows from \eqref{pierewasze szacownaiw}, \eqref{normalisation1} and \eqref{warunek na miary radialne} that
\begin{align}
\int_{\mathbb{V}_{\eta}}\int_{0}^{+\ns}H\rbr{b\cdot r\left\langle G(x),\xi\right\rangle }\gamma_\xi\rbr{\dd r}\lambda(\dd \xi)
 & \le4\int_{\mathbb{V}_{\eta}}\int_{0}^{+\ns}H\rbr{b\cdot r\left|G(x)\right|\eta}\gamma_\xi\rbr{\dd r}\lambda(\dd \xi)\nonumber \\
 & \le4\varepsilon\sup_{\xi\in\mathbb{V}_{\eta}}\int_{0}^{+\ns}H\rbr{b\cdot r\left|G(x)\right|\eta}\gamma_\xi\rbr{\dd r}\nonumber \\
 & \le4 K \varepsilon\inf_{\xi\in\mathbb{S}^{d-1}\setminus\mathbb{V}_{\eta}}\int_{0}^{+\ns}H\rbr{b\cdot r\left|G(x)\right|\eta}\gamma_\xi\rbr{\dd r}.\label{eq:osz_kur}
\end{align}
From Lemma \ref{lem oszacownaia na funkcje podcalkowa}, for $b,r\ge0$ and $\xi\in\S^{d-1}$ such that $\left\langle G_{\ns},\xi\right\rangle \in[\eta,1]$, 
we also estimate 
\begin{align}
 H\rbr{b\cdot r\left\langle G(x),\xi\right\rangle } & \le H\rbr{b\cdot r\left|G(x)\right|\left(\left\langle G_{\ns},\xi\right\rangle +\left|\left\langle \frac{G(x)}{\left|G(x)\right|}-G_{\ns},\xi\right\rangle \right|\right)}\nonumber \\
 & \le H\rbr{b\cdot r\left|G(x)\right|\left(\left\langle G_{\ns},\xi\right\rangle +\left\langle G_{\ns},\xi\right\rangle \varepsilon\right)}\nonumber \\
 & \le\left(1+\varepsilon\right)^{2}H\rbr{b\cdot r\left|G(x)\right|\left\langle G_{\ns},\xi\right\rangle },\label{eq:osz_kur_3}
\end{align}
and 
\begin{align}
H\rbr{b\cdot r\left\langle G(x),\xi\right\rangle }\nonumber 
 & \ge H\rbr{b\cdot r\left|G(x)\right|\left(\left\langle G_{\ns},\xi\right\rangle -\left|\left\langle \frac{G(x)}{\left|G(x)\right|}-G_{\ns},\xi\right\rangle \right|\right)}\nonumber \\
 & \ge H\rbr{b\cdot r\left|G(x)\right|\left(\left\langle G_{\ns},\xi\right\rangle -\left\langle G_{\ns},\xi\right\rangle \varepsilon\right)}\nonumber \\
 & \ge\left(1-\varepsilon\right)^{2}H\rbr{b\cdot r\left|G(x)\right|\left\langle G_{\ns},\xi\right\rangle }.\label{eq:osz_kur_2}
\end{align}
Notice that by \eqref{normalisation} and \eqref{normalisation1}, $\lambda\rbr{\mathbb{S}^{d-1} \setminus \mathbb{V}_{\eta}} \ge 1-\varepsilon$. From \eqref{eq:osz_kur_2} and then from  $\lambda\rbr{\mathbb{S}^{d-1}\setminus \mathbb{V}_{\eta}} \ge 1-\varepsilon$ and 
 \eqref{eq:osz_kur} we obtain 
\begin{align}\nonumber
&\int_{\S^{d-1}\setminus\mathbb{V}_{\eta}}\int_{0}^{+\ns}H\rbr{b\cdot r\left\langle G(x),\xi\right\rangle }\gamma_\xi\rbr{\dd r}\lambda(\dd \xi)\\ \nonumber
 & \ge\int_{\S^{d-1}\setminus\mathbb{V}_{\eta}}\int_{0}^{+\ns}\left(1-\varepsilon\right)^{2}H\rbr{b\cdot r\left|G(x)\right|\left\langle G_{\ns},\xi\right\rangle }\gamma_\xi\rbr{\dd r}\lambda(\dd \xi)\nonumber \\
 & \ge\left(1-\varepsilon\right)^{2}\int_{\S^{d-1}\setminus\mathbb{V}_{\eta}}\int_{0}^{+\ns}H\rbr{b\cdot r\left|G(x)\right|\eta}\gamma_\xi\rbr{\dd r}\lambda(\dd \xi)\nonumber \\
 & \ge\left(1-\varepsilon\right)^{2}\left(1-\varepsilon\right)\inf_{\xi\in\S^{d-1}\setminus\mathbb{V}_{\eta} }\int_{0}^{+\ns}H\rbr{b\cdot r\left|G(x)\right|\eta}\gamma_{\xi}\rbr{\dd r}\nonumber \\
 & \ge\frac{\left(1-\varepsilon\right)^{3}}{4 K \varepsilon}\int_{\mathbb{V}_{\eta}}\int_{0}^{+\ns}H\rbr{b\cdot r\left\langle G(x),\xi\right\rangle }\gamma_{\xi}\rbr{\dd r}\lambda(\dd \xi).\label{eq:osz_kur_4}
\end{align}
From (\ref{eq:osz_kur_4}) and (\ref{eq:osz_kur_3}) we estimate
\begin{align*}
J_{X}\rbr{bG(x)}= & \int_{\S^{d-1}\setminus\mathbb{V}_{\eta}}\int_{0}^{+\ns}H\rbr{b\cdot r\left\langle G(x),\xi\right\rangle }\gamma_{\xi}\rbr{\dd r}\lambda(\dd \xi)+\int_{\mathbb{V}_{\eta}}\int_{0}^{+\ns}H\rbr{b\cdot r\left\langle G(x),\xi\right\rangle }\gamma_{\xi}\rbr{\dd r}\lambda(\dd \xi)\\
\le & \int_{\S^{d-1}\setminus\mathbb{V}_{\eta}}\int_{0}^{+\ns}H\rbr{b\cdot r\left\langle G(x),\xi\right\rangle }\gamma_{\xi}\rbr{\dd r}\lambda(\dd \xi)\\
 & +\frac{4 K \varepsilon}{\left(1-\varepsilon\right)^{3}}\int_{\S^{d-1}\setminus\mathbb{V}_{\eta}}\int_{0}^{+\ns}H\rbr{b\cdot r\left\langle G(x),\xi\right\rangle }\gamma_{\xi}\rbr{\dd r}\lambda(\dd \xi)\\
\le & \left(1+\varepsilon\right)^{2}\left(1+\frac{4 K \varepsilon}{\left(1-\varepsilon\right)^{3}}\right)\int_{\S^{d-1}}\int_{0}^{+\ns}H\rbr{b\cdot r\left|G(x)\right|\left\langle G_{\ns},\xi\right\rangle }\gamma_{\xi}\rbr{\dd r}\lambda(\dd \xi)\\
= & \left(1+\varepsilon\right)^{2}\left(1+\frac{4 K \varepsilon}{\left(1-\varepsilon\right)^{3}}\right)M_{G_{\ns}}\left(b\cdot r\left|G(x)\right|\right).
\end{align*}
Hence 
\begin{equation}
J_{X}\rbr{bG(x)}\le\left(1+\varepsilon\right)^{2}\left(1+\frac{4 K \varepsilon}{\left(1-\varepsilon\right)^{3}}\right)M_{G_{\ns}}\left(b\cdot r\left|G(x)\right|\right).\label{eq:mniam}
\end{equation}
In order to get the lower bound let us notice that 
\begin{align*}
\int_{\S^{d-1}\setminus\mathbb{V}_{\eta}}\int_{0}^{+\ns}H\rbr{b\cdot r\left|G(x)\right|\left\langle G_{\ns},\xi\right\rangle}\gamma_{\xi}\rbr{\dd r}\lambda(\dd \xi)
 & \ge\int_{\S^{d-1}\setminus\mathbb{V}_{\eta}}\int_{0}^{+\ns}H\rbr{b\cdot r\left|G(x)\right|\eta}\gamma_{\xi}\rbr{\dd r}\lambda(\dd \xi)\\
 & \ge\left(1-\varepsilon\right)\inf_{\xi\in\S^{d-1}\setminus\mathbb{V}_{\eta}}\int_{0}^{+\ns}H\rbr{b\cdot r\left|G(x)\right|\eta}\gamma_{\xi}\rbr{\dd r},
\end{align*}
and 
\begin{align*}
 \int_{\mathbb{V}_{\eta}}\int_{0}^{+\ns}H\rbr{b\cdot r\left|G(x)\right|\left\langle G_{\ns},\xi\right\rangle }\gamma_{\xi}\rbr{\dd r}\lambda(\dd \xi) & \le\int_{\mathbb{V}_{\eta}}\int_{0}^{+\ns}H\rbr{b\cdot r\left|G(x)\right|\eta}\gamma_{\xi}\rbr{\dd r}\lambda(\dd \xi)\\
 & \le\varepsilon\sup_{\xi\in\mathbb{V}_{\eta}}\int_{0}^{+\ns}H\rbr{b\cdot r\left|G(x)\right|\eta}\gamma_{\xi}\rbr{\dd r}\\
  & \le K \varepsilon\inf_{\xi\in\S^{d-1}\setminus\mathbb{V}_{\eta}}\int_{0}^{+\ns}H\rbr{b\cdot r\left|G(x)\right|\eta}\gamma_{\xi}\rbr{\dd r}.
\end{align*}
Hence
\begin{align}
 & \int_{\S^{d-1}\setminus\mathbb{V}_{\eta}}\int_{0}^{+\ns}H\rbr{b\cdot r\left|G(x)\right|\left\langle G_{\ns},\xi\right\rangle }\gamma_{\xi}\rbr{\dd r}\lambda(\dd \xi)\nonumber \\
 & \ge\frac{1-\varepsilon}{K \varepsilon}\int_{\mathbb{V}_{\eta}}\int_{0}^{+\ns}H\rbr{b\cdot r\left|G(x)\right|\left\langle G_{\ns},\xi\right\rangle }\gamma_{\xi}\rbr{\dd r}\lambda(\dd \xi),\label{eq:osz_kur_5}
\end{align}
and from this we obtain 
\begin{align}
\int_{\S^{d-1}}\int_{0}^{+\ns}&H\rbr{b\cdot r\left|G(x)\right|\left\langle G_{\ns},\xi\right\rangle }\gamma_{\xi}\rbr{\dd r}\lambda(\dd \xi) =\int_{\S^{d-1}\setminus\mathbb{V}_{\eta}}\int_{0}^{+\ns}H\rbr{b\cdot r\left|G(x)\right|\left\langle G_{\ns},\xi\right\rangle }\gamma_{\xi}\rbr{\dd r}\lambda(\dd \xi) \nonumber \\
 & \quad+\int_{\mathbb{V}_{\eta}}\int_{0}^{+\ns}H\rbr{b\cdot r\left|G(x)\right|\left\langle G_{\ns},\xi\right\rangle }\gamma_{\xi}\rbr{\dd r}\lambda(\dd \xi)\label{eq:osz_kur_5-1} \nonumber \\
 & \le\left(1+\frac{K \varepsilon}{1-\varepsilon}\right)\int_{\S^{d-1}\setminus\mathbb{V}_{\eta}}\int_{0}^{+\ns}H\rbr{b\cdot r\left|G(x)\right|\left\langle G_{\ns},\xi\right\rangle }\gamma_{\xi}\rbr{\dd r}\lambda(\dd \xi).
\end{align}
From (\ref{eq:osz_kur_2}) and (\ref{eq:osz_kur_5-1}) we
get 
\begin{align*}
J_{X}\rbr{bG(x)}\ge & \int_{\S^{d-1}\setminus\mathbb{V}_{\eta}}\int_{0}^{+\ns}H\rbr{b\cdot r\left\langle G(x),\xi\right\rangle }\gamma_{\xi}\rbr{\dd r}\lambda(\dd \xi)\\
\ge & \left(1-\varepsilon\right)^{2}\int_{\S^{d-1}\setminus\mathbb{V}_{\eta}}\int_{0}^{+\ns}H\rbr{b\cdot r\left|G(x)\right|\left\langle G_{\ns},\xi\right\rangle }\gamma_{\xi}\rbr{\dd r}\lambda(\dd \xi)\\
\ge & \frac{\left(1-\varepsilon\right)^{2}}{\left(1+\frac{K \varepsilon}{1-\varepsilon}\right)}\int_{\S^{d-1}}\int_{0}^{+\ns}H\rbr{b\cdot r\left|G(x)\right|\left\langle G_{\ns},\xi\right\rangle }\gamma_{\xi}\rbr{\dd r}\lambda(\dd \xi)\\
= & \frac{\left(1-\varepsilon\right)^{2}}{\left(1+\frac{K \varepsilon}{1-\varepsilon}\right)}M_{G_{\ns}}\left(b\cdot r\left|G(x)\right|\right).
\end{align*}
Hence 
\begin{equation}
J_{X}\rbr{bG(x)}\ge\frac{\left(1-\varepsilon\right)^{2}}{\left(1+\frac{ K \varepsilon}{1-\varepsilon}\right)}M_{G_{\ns}}\left(b\cdot r\left|G(x)\right|\right).\label{eq:mniam_mniam}
\end{equation}
Now (\ref{eq:teeza}) follows from (\ref{eq:mniam}), (\ref{eq:mniam_mniam})
and (\ref{eq:mniam_mniam_mniam}).
\hfill$\square$

\begin{lem}\label{lem o charakteryzacji f potegowej}
Let $J_{\rho}$ be given by \eqref{funkcja Jot} with $\rho(\dd v)$ satisfying \eqref{Ro}. Assume that 
\begin{gather}\label{beta war do stabilnosci}
J_\rho(\beta b)=\eta J_{\rho}(b), \quad b\geq 0,
\end{gather}
and 
\begin{gather}\label{gamma war do stabilnosci}
J_\rho(\gamma b)=\theta J_{\rho}(b), \quad b\geq 0,
\end{gather}
for some $\beta>1$, $\gamma>1$ such that $\ln\beta/\ln\gamma\notin \mathbb{Q}$ and $\eta>1$, $\theta>1$. Then
\begin{gather}\label{wzor do udowodnienia}
J_\rho(b)=C b^{\alpha}, b\geq 0,
\end{gather}
for some $C>0$ and $\alpha\in(1,2)$.
\end{lem}
{\bf Proof :} By iterative application of \eqref{beta war do stabilnosci} and \eqref{gamma war do stabilnosci} we see that for any 
$m,n\in\mathbb{N}$
\begin{gather*}
J_{\rho}(\beta^m\gamma^n b)=\eta^m\theta^n J_\rho(b), \quad b\geq 0,
\end{gather*}
which can be written as
\begin{gather}\label{wyjsciowa rownosc w dowodzie}
J_{\rho}(b e^{m \ln\beta+n \ln\gamma})=e^{m\ln\eta+n \ln\theta} J_\rho(b), \quad b\geq 0.
\end{gather}
In Lemma \ref{lem Weyl} below we prove that the set 
$$
D:=\{m \ln\beta -n\ln\gamma; \ m,n\in\mathbb{Z}\}
$$
is dense in $\mathbb{R}$. So, for any $\delta>0$ there exist $m,n\in\mathbb{Z}, m\neq 0,$ such that
\begin{gather}\label{odleglosc 1}
\mid m \ln\beta -n\ln\gamma\mid<\delta,
\end{gather}
and then, by \eqref{min maxxx} and \eqref{wyjsciowa rownosc w dowodzie}, we obtain that
\begin{gather}\label{odleglosc 2}
e^{-2 \delta}\leq \frac{e^{m\ln \eta}}{e^{n\ln \theta}}=\frac{J_\rho(e^{m\ln \beta})}{J_\rho(e^{n\ln \gamma})}\leq e^{2\delta}.
\end{gather}
It follows from \eqref{odleglosc 1} that
$$
\left\vert\frac{\ln\beta}{\ln\gamma}-\frac{n}{m}\right\vert\leq \frac{\delta}{\mid m\mid\ln\gamma},
$$
and from \eqref{odleglosc 2} that
$$
\left\vert\frac{\ln\eta}{\ln\theta}-\frac{n}{m}\right\vert\leq \frac{2\delta}{\mid m\mid\ln\theta}.
$$
Consequently,
$$
\left\vert\frac{\ln\beta}{\ln\gamma}-\frac{\ln\eta}{\ln\theta}\right\vert\leq\frac{\delta}{\mid m\mid\ln\gamma}+
\frac{2\delta}{\mid m\mid\ln\theta}\leq\frac{\delta}{\ln\gamma}+
\frac{2\delta}{\ln\theta}.
$$
Letting $\delta\longrightarrow 0$ yields
$$
\frac{\ln\beta}{\ln\gamma}=\frac{\ln\eta}{\ln\theta}.
$$
Let us define
$$
\alpha:=\frac{\ln\eta}{\ln \beta}=\frac{\ln\theta}{\ln\gamma}>0,
$$
and put $b=1$ in \eqref{wyjsciowa rownosc w dowodzie}. This gives
$$
J_\rho(e^{m \ln\beta+n \ln\gamma})=J_\rho(1)\left(e^{m \ln\beta+n \ln\gamma}\right)^\alpha,
$$
which means that $J_\rho(b)=J_\rho(1)b^\alpha$ for $b$ from the set $e^D$
which is dense in $[0,+\infty)$. As $J_\rho$ is continuous, \eqref{wzor do udowodnienia} follows.
Finally, by Proposition 3.4 in \cite{BarskiLochowski} it follows that $\alpha\in(1,2)$. 
\hfill$\square$

\vskip1ex
The following result is strictly related to Weyl's equidistribution theorem, see \cite{Weyl}.

\begin{lem}\label{lem Weyl}
Let $p,q >0$ be such that $p/q\notin\mathbb{Q}$. Let us define the set
$$
G:=\{mp+nq; \quad m,n,=1,2,...\}.
$$
Then for each $\delta>0$ there exists a number $M(\delta)>0$ such that
$$
\forall x\geq M(\delta) \quad  \exists \ g\in G \quad\text{such that} \ \mid x-g\mid\leq \delta.
$$
Moreover, the set
$$
D:=\{mp+nq; \quad m,n\in\mathbb{Z}\},
$$
is dense in $\mathbb{R}$.
\end{lem}
{\bf Proof:} Since $p/q\notin\mathbb{Q}$, at least one of $p,q$, say $q$, is irrational. For simplicity assume that $p=1$ and consider the sequence 
$$
r(jq), j=1,2,... \quad \text{where} \  r(x):=x \ \text{mod} \ 1,
$$
of fractional parts of the numbers $jq, j=1,2,...$\ . Recall, Weyl's equidistribution theorem states that 
\begin{gather}\label{weyl result}
\lim_{N\longrightarrow +\infty}\frac{\sharp\{j\leq N: r(jq)\in[a,b]\}}{N}=b-a
\end{gather}
for any $[a,b]\subseteq [0,1)$ if and only if $q$ is irrational.

For fixed $\delta>0$ and $n$ such that $1/n<\delta$ let us consider a partition of $[0,1)$ of the form
$$
[0,1)=\bigcup_{k=0}^{n-1}A_k, \quad A_k:=[k/n, (k+1)/n).
$$
For a natural number $N$ let us consider the set $R_N:=\{r(jq) : j=1,2,...,N\}$. By \eqref{weyl result}, for each $k=0,1,...,n-1$, there exists $N_k$ such that 
$$
R_{N_k}\cap A_k\neq \emptyset.
$$
Then for $\bar{N}:=\max\{N_0,N_1,...,N_{n-1}\}$ we have
$$
R_{\bar{N}}\cap A_k\neq \emptyset, \quad k=0,1,...,n-1.
$$
Let $M=M(\delta):=\bar{N}q$. Then, for $x\geq M$, there exists a number $N_x\leq \bar{N}$ such that
\begin{gather}\label{oszcownaie reszt}
\mid r(N_x q)-r(x)\mid\leq \frac{1}{n}.
\end{gather}
Then for the number 
$$
g:=\lfloor x \rfloor-\lfloor N_x q\rfloor+N_xg\in G
$$
the following holds
\begin{align*}
\mid x-g\mid&=\mid x-(\lfloor x \rfloor-\lfloor N_x q\rfloor+N_xq)\mid\\[1ex]
&=\mid\lfloor x \rfloor+r(x)-\lfloor x \rfloor+\lfloor N_x q\rfloor-N_xq\mid\\[1ex]
&=\mid r(x)-r(N_xq)\mid\leq 1/n<\delta,
\end{align*}
where the last inequality follows from \eqref{oszcownaie reszt}.

The density of $D$ is an immediate consequence of the first part of the Lemma. Indeed, for $x<M(\delta)$ and $g\in G$ such that $x+g>M(\delta)$ there exists $\tilde{g}\in G$ such that $\mid x+g-\tilde{g}\mid<\delta$.

The general case with $p\neq 1$ can be proven in the same way but requires a generalized version of Weyl's theorem, which says that the numbers $r_p(nq), n=1,2,...$, where $r_p(x):=x \ mod \ p$, are equidistributed on $[0,p)$ if and only if $p/q\notin \mathbb{Q}$. This can be proven by a straightforward modification of the original arguments of Weyl.
\hfill $\square$

\subsection{Proof of Theorem \ref{tw_sferyczne}}  \label{proof_of_tw_sf}

By Remark \ref{rem o rownaniu dla Laplasuf} and Remark \ref{rem o charakterystykach dla alfa stab. 1 dim} the Laplace transform $J_X$ satisfies
\begin{equation} \label{warrr1}
J_{X}(b G(x))=J_{\nu_{G(0)}}(b)+xJ_{\mu}(b), \quad b,x\geq 0,
\end{equation}
where $\mu(\dd v)$ is the measure satisfying \eqref{war calkowe na mu}-\eqref{rozklad nu G(x)}. By discussion preceding the formulation of Theorem \ref{tw_sferyczne}  we have $G(0) = 0$, hence  \eqref{warrr1} simpilfes to
\begin{equation} \label{warrr2}
J_{X}(b G(x))=xJ_{\mu}(b), \quad x\geq 0.
\end{equation}

Assumption \eqref{spherical part non-degenerated} and \eqref{warrr2} imply that $J_{X}(y), J_{\mu}(b)>0$, $G(x) \neq 0$, for $y \in \R^d \setminus \cbr{0}, b>0$, $x>0$.

Let $G_{0} = \lim_{x \ra 0+} \frac{G(x)}{|G(x)|}$. 
It follows from Proposition \ref{sferycznosc} that there exists a function $\delta:(0,+\infty)\rightarrow(0,+\infty)$,
such that for any $\varepsilon>0$ from the inequality 
\[
\left|\frac{G(x)}{\left|G(x)\right|}-G_{0}\right|\le\delta(\varepsilon).
\]
follows that for any $b\ge0$
\[
1-\varepsilon\le\frac{J_{X}\rbr{b\frac{G(x)}{\left|G(x)\right|}}}{J_{X}\rbr{bG_{0}}}\le1+\varepsilon.
\]
Thus for any $\varepsilon>0$ there exists $m(\varepsilon) >0$,
such that for $x\in \rbr{0, m(\varepsilon)}$
\[
\left|\frac{G(x)}{\left|G(x)\right|}-G_{0}\right|\le\delta(\varepsilon),
\]
and hence for any $b>0$
\[
1-\varepsilon\le\frac{J_{X}\rbr{b\frac{G(x)}{\left|G(x)\right|}}}{J_{X}\rbr{bG_{0}}}\le1+\varepsilon.
\]
Let us fix $\beta>1$ and take $x_{1},x_{2}$ satisfying
$0 < x_{1}\le x_{2} < m(\varepsilon)$, $\beta\left|G(x_{1})\right|=\left|G(x_{2})\right| > 0$
(from the continuity of $G$  it follows that such $x_{1}$ and $x_{2}$ exist).
Then for any $b>0$ and $i=1,2$, by \eqref{warrr2}, 
\[
1-\varepsilon\le\frac{J_{X}\rbr{b\frac{G(x_{i})}{\left|G(x_{i})\right|}}}{J_{X}\rbr{bG_{0}}}=\frac{x_{i}J_{\mu}\rbr{\frac{b}{\left|G(x_{i})\right|}}}{J_{X}\rbr{bG_{0}}}\le1+\varepsilon.
\]
Hence for any $b>0$, taking $\tilde{b}=\beta\left|G(x_{1})\right|b$ we get
\[
\frac{1-\varepsilon}{1+\varepsilon}\cdot \frac{x_{2}}{x_{1}}\le\frac{J_{\mu}\rbr{\frac{\tilde{b}}{\left|G(x_{1})\right|}}}{J_{\mu}\rbr{\frac{\tilde{b}}{\left|G(x_{2})\right|}}}=\frac{J_{\mu}\rbr{\beta b}}{J_{\mu}\rbr b}\le\frac{1+\varepsilon}{1-\varepsilon}\cdot \frac{x_{2}}{x_{1}}
\]
which yields
\[
\frac{1-\varepsilon}{1+\varepsilon}\cdot \frac{J_{\mu}\rbr{\beta b}}{J_{\mu}\rbr b}\le\frac{x_{2}}{x_{1}}\le\frac{1+\varepsilon}{1-\varepsilon}\cdot \frac{J_{\mu}\rbr{\beta b}}{J_{\mu}\rbr b}.
\]
Since $\varepsilon>0$ is arbitrary, taking $\varepsilon\rightarrow0$
and $x_{1},x_{2}$ satisfying
$0 < x_{1}\le x_{2} <m(\varepsilon)$, $\beta\left|G(x_{1})\right|=\left|G(x_{2})\right|$ we obtain that 
\[
\lim_{\varepsilon\rightarrow0}\frac{x_{2}}{x_{1}}=\eta,
\]
where $\eta = J_{\mu}\rbr{\beta b}/J_{\mu}\rbr b > 1$ is independent of $b>0$.
Hence, for all $b\ge0$ we have 
\[
J_{\mu}\rbr{\beta b}=\eta J_{\mu}\rbr b.
\]

Similarly, take $\gamma>1$ such that $\ln\beta/\ln\gamma\notin \mathbb{Q}$. Reasoning similarly as before we get that there exists $\theta>1$, such that for all $b\ge0$ we have
\[
J_{\mu}\rbr{\gamma b}=\theta J_{\mu}\rbr b.
\]
Now the thesis follows from Lemma \ref{lem o charakteryzacji f potegowej} and the one to one correspondence
between Laplace transforms and measures on $[0,+\ns)$, see \cite{{Feller}} p. 233.
\hfill $\square$

\end{document}